\documentclass[12pt,a4paper,leqno]{article}

\usepackage{latexsym}
\usepackage{amssymb,exscale}
\usepackage[centertags]{amsmath}
\usepackage{amsthm}
\usepackage[all]{xy}

\numberwithin{equation}{section}

\swapnumbers
\theoremstyle{definition}
\newtheorem{theorem}[equation]{Theorem}
\newtheorem{lemma}[equation]{Lemma}
\newtheorem{proposition}[equation]{Proposition}
\newtheorem{corollary}[equation]{Corollary}
\newtheorem{conjecture}[equation]{Conjecture}
\newtheorem{definition}[equation]{Definition}

\newtheorem{note}[equation]{Note}
\newtheorem{myItem}[equation]{}

\renewcommand{\phi}{\varphi}

\renewcommand{\(}{\bigl(}
\renewcommand{\)}{\bigr)\vphantom{)}}

\newcommand{\Prob}{\operatorname{Prob}}
\newcommand{\Homeo}{\operatorname{Homeo}}
\newcommand{\entier}{\operatorname{entier}}
\newcommand{\Aut}{\operatorname{Aut}}
\newcommand{\Iso}{\operatorname{Iso}}
\newcommand{\MALG}{\operatorname{MALG}}
\newcommand{\dist}{\operatorname{dist}}

\newcommand{\SF}{{\text{S$\sigma$F}}}
\newcommand{\K}{\mathrm K}
\newcommand{\bF}{\mathrm F}
\newcommand{\bL}{\mathrm L}
\newcommand{\U}{\mathrm U}

\newcommand{\Om}{\Omega}
\newcommand{\om}{\omega}
\newcommand{\al}{\alpha}
\newcommand{\be}{\beta}
\newcommand{\eps}{\varepsilon}
\newcommand{\ga}{\gamma}
\newcommand{\si}{\sigma}
\newcommand{\de}{\delta}

\newcommand{\F}{\mathcal F}
\newcommand{\E}{\mathcal E}

\newcommand{\M}{\mathcal M}
\newcommand{\fav}{{\text{fav}}}
\newcommand{\Omfav}{\Om^\fav}
\newcommand{\Ffav}{\F^\fav}
\newcommand{\Pfav}{P^\fav}
\newcommand{\Mfav}{\M^\fav}
\newcommand{\Gfav}{G^\fav}
\newcommand{\Xcna}{X^{\text{cna}}}

\newcommand{\Ex}{\mathbb E}
\renewcommand{\Pr}[1]{\,\mathbb P\,\(\,#1\,\)\,}
\newcommand{\R}{\mathbb R}

\newcommand{\Q}{\mathbb Q}
\newcommand{\cE}[2]{\mathbb{E}\,\(\,#1\,\big|\,#2\,\)\,}
\newcommand{\scE}[2]{\,\mathbb{E}(#1|#2)}

\newcommand{\cP}[2]{\mathbb{P}\,\(\,#1\,\big|\,#2\,\)\,}

\newcommand{\ti}{\tilde}

\newcommand{\sif}{$\sigma$-field}

\newcommand{\Zero}{\mathbf0}
\newcommand{\One}{\mathbf1}
\newcommand{\imply}{\;\;\;\Longrightarrow\;\;\;}

\newcommand{\imp}{$ \Longrightarrow $ }
\newcommand{\impl}{\;\Longrightarrow\;}
\renewcommand{\equiv}{\;\;\;\Longleftrightarrow\;\;\;}

\def\emailwww#1#2{\par\qquad {\tt #1}\par\qquad {\tt #2}\medskip}

\newenvironment{myitemize}{\begin{list}{$\bullet$}
{\setlength{\topsep}{1mm}
\setlength{\partopsep}{0mm}
\setlength{\itemsep}{1mm}
\setlength{\parsep}{0mm}
\setlength{\parskip}{0mm}}}
{\end{list}}

\begin{document}

\title{Filtrations of random processes\\
  in the light of classification theory.\\
  I. A topological zero-one law.}
\author{Boris Tsirelson}
\date{}

\maketitle

\stepcounter{footnote}
\footnotetext{%
 Supported by the \textsc{israel science foundation} founded by the
 Israel Academy of Sciences and Humanities.}

\begin{abstract}
Filtered probability spaces (called ``filtrations'' for short) are
shown to satisfy such a topological zero-one law: for every property
of filtrations, either the property holds for almost all filtrations,
or its negation does. In particular, almost all filtrations are
conditionally nonatomic.

An accurate formulation is given in terms of orbit equivalence
relations on Polish $ G $-spaces. The set of all isomorphic classes of
filtrations may be identified with the orbit space $ X/G $ (such a
space is sometimes called a \emph{singular space} in the modern
descriptive set theory, see Kechris \cite[\S 2]{Ke99}) for a special
Polish $ G $-space $ X $. A ``property of filtrations'' means a $ G
$-invariant subset of $ X $ having the Baire property. ``Almost all
filtrations'' means a comeager subset of $ X $ (the Baire category
approach). The zero-one law is a kind of ergodicity of $ X $. It holds
for filtrations both in discrete and continuous time.

The interplay between probability theory and descriptive set theory
could be interesting for both parties.
\end{abstract}

\section{Preliminaries on the classification theory (for probabilists)}
\subsection{Topological counterparts of some probabilistic ideas}

Objects to be classified (say, filtrations) often form an
infinite-dimensional space with a natural topology but without a
natural measure. We cannot say `for almost all filtrations' in the
measure-theoretic sense. Instead we can say `for comeager many
filtrations' in the sense explained below. Fubini theorem, zero-one
law and other nice probabilistic tools have topological counterparts
based on the \emph{Baire category theorem.}

\begin{definition}
(a) A \emph{Polish metric space} is a complete separable metric
space.

(b) A \emph{Polish topological space} (or just `Polish space') is a
topological space that admits a compatible metric turning it into a
Polish metric space.
\end{definition}

(`Compatible' means that the topology corresponding to the metric is
the same as the given topology.)

\begin{theorem}(Baire; see \cite[2.5.5]{Sri} or \cite[8.4 and
8.1(iii)]{Ke} or \cite[2.5.2]{Du})
Let $ U_1, U_2, \dots $ be dense open sets in a Polish space $ X
$. Then their intersection $ U_1 \cap U_2 \cap \dots $ is dense in $ X
$. 
\end{theorem}

\begin{definition}\label{1.3}
Let $ X $ be a Polish space and $ A \subset X $ a set.

(a) $ A $ is \emph{comeager,} if $ A \supset U_1 \cap U_2 \cap \dots $
for some dense open sets $ U_1, U_2, \dots $

(b) $ A $ is \emph{meager} (or `of first category') if $ X \setminus A
$ is comeager, that is, if $ A \subset F_1 \cup F_2 \cup \dots $ for
some closed sets $ F_1, F_2, \dots $ with no interior points.

(c) `For \emph{comeager many} $ x \in X $' (symbolically, $ \forall^*
x \in X $) means: for all $ x $ of a comeager set $ A \subset X $.

(d) `For \emph{non-meager many} $ x \in X $' (symbolically, $ \exists^*
x \in X $) means: for all $ x $ of a non-meager set $ A \subset X $.
\end{definition}

(See \cite[Sect.~8.J]{Ke}.)
The Baire category theorem ensures that a comeager set cannot be
meager. Thus,
\[
\forall x \, (\dots) \imply \forall^* x \, (\dots) \imply \exists^* x
\, (\dots) \imply \exists x \, (\dots) \, .
\]

We know that Lebesgue measurable sets are Borel sets modulo sets of
measure zero. Similarly, sets having the Baire property are Borel sets
modulo meager sets, which is equivalent to the next definition.

\begin{definition}
Let $ X $ be a Polish space.

(a) A set $ A \subset X $ has the \emph{Baire property} if there
exists an open set $ U \subset X $ such that $ A \setminus U $ and $ U
\setminus A $ are meager.

(b) A function $ f $ from $ X $ to another Polish space $ Z $ is
called \emph{Baire measurable} if $ f^{-1} (V) $ has the Baire
property for every open set $ V \subset Z $.
\end{definition}

\begin{proposition}
Sets having the Baire property are a \sif.\footnote{%
 In other words, $ \si $-algebra.}
It is the \sif\ generated by all open (or equivalently, Borel) sets
 and all meager sets. (See \cite[8.22]{Ke}.)\footnote{%
 Baire property holds for all Borel sets, and moreover, for all
 analytic sets (Lusin-Sierpinski theorem, see \cite[21.6]{Ke}). For
 projective (and even more general) sets, Baire property is proven
 using additional (to ZFC) axioms (see \cite[26.3, 36.20; 38.15,
 38.17(ii); 39.17]{Ke}), which is considered a fault
 of ZFC axiomatics (see \cite[p.~16]{Fo}), not of Baire property!}
\end{proposition}

A surprise: Borel sets are open modulo meager. In contrast, they are $
G_\de $ (or $ F_\si $, but generally not open) sets modulo sets of
measure zero. The boundary of an open set is always meager, but not
always of measure zero.

Here is a counterpart of Fubini theorem, and then --- of Kolmogorov's
zero-one law.

\begin{theorem}(Kuratowski-Ulam; see \cite[8.41(iii)]{Ke} or
\cite[3.5.16]{Sri})
Let $ X,Y $ be Polish spaces, and a set $ A \subset X \times Y $ have
the Baire property. Then\footnote{%
 Saying ``$ \forall^* x \, \forall^* y $'' we mean that the comeager
 set of these $ y $ may depend on $ x $.}
\[
\forall^* x \, \forall^* y \; (x,y) \in A \equiv A \text{ is comeager}
\equiv \forall^* y \, \forall^* x \; (x,y) \in A \, .
\]
\end{theorem}

\begin{theorem}(A topological zero-one law; see \cite[8.47]{Ke}.)
Let $ X_1, X_2, \dots $ be Polish spaces, and a set $ A \subset X_1
\times X_2 \times \dots $ have the Baire property. If $ A $ is a tail
set\footnote{%
 That is, the relation $ (x_1,x_2,\dots) \in A $ is insensitive to any
 change of a finite number of $ x_k $.}
then $ A $ is either meager or comeager.
\end{theorem}

\subsection{Constructing Polish spaces}

\begin{myItem}
Let $ (\Om,\F,P) $ be a probability space, $ L_2 (\Om,\F,P) $ the
Hilbert space of square integrable random variables,\footnote{%
 By a random variable I mean here an equivalence class of measurable
 functions $ \Om \to \R $.}
$ L_0 (\Om,\F,P)
$ the linear topological space of all random variables (equipped with
the topology of convergence in probability), $ \MALG (\Om,\F,P) $ the
set (in fact, complete Boolean algebra) of equivalence classes of $ \F
$-measurable sets, equipped with the metric $ \dist (A,B) = P ( A
\setminus B ) + P ( B \setminus A ) $. Assume that $ (\Om.\F,P) $ is
\emph{separable} in the sense that one of the three spaces ($ L_2 $, $ L_0 $,
$ \MALG $) is separable, then others are also separable. All the three
spaces are Polish ($ L_2 $ and $ \MALG $ being Polish \emph{metric}
spaces).\footnote{%
 Also $ L_0 $ becomes a Polish \emph{metric} space, being equipped
 with one of several well-known equivalent metrics, see
 \cite[17.45]{Ke}, \cite[9.2.3]{Du}.}
\end{myItem}

\begin{myItem}\label{1.9a}
Let $ (\Om,\F,P) $ be a separable probability space and $ X $ a Polish
space, then $ X $-valued random variables on $ (\Om,\F,P) $ form a
Polish space $ L_0 (\Om,\F,P; \linebreak[0] X) $ (but not a linear space, in
general). Equivalent metrics on $ X $ induce equivalent metrics on $
L_0 (\Om,\F,P; X) $, which follows easily from the well-known relation
between convergence in probability and convergence almost everywhere;
namely, $ f_n \to f $ in probability if and only if every subsequence
of $ (f_n) $ contains a subsequence converging to $ f $ almost
everywhere. See \cite[9.2.1]{Du}.
\end{myItem}

\begin{myItem}\label{1.10a}
Let $ X $ be a Polish space. Denote by $ \Prob (X) $ the set of
all probability distributions on $ X $ (that is, probability measures
on the Borel \sif\ of $ X $), equipped with the weak
topology, called also the weak$^*$ topology. It is the unique metrizable
topology such that for every $ \mu, \mu_1, \mu_2, \dots \in \Prob (X) $
\[
\mu_n \to \mu \quad \text{if and only if} \quad \forall \phi \;\;
 \int \phi \, d\mu_n \to \int \phi \, d\mu \, ,
\]
where $ \phi $ runs over all bounded continuous functions $ X \to \R
$. See also \cite[Sect.~9.3]{Du}, \cite[17.20]{Ke}. Some countable set of
functions $ \phi $ is enough for generating the same topology, see
\cite[17.20(ii)]{Ke}.

The topological space $ \Prob (X) $ is Polish. See
\cite[17.23]{Ke}. Note also the continuous map
\[
L_0 (\Om,\F,P; X) \ni f \mapsto P_f \in \Prob(X) \, ;
\]
the distribution $ P_f $ of $ f $ is defined by $ P_f (A) = P (
f^{-1}(A) ) $. See \cite[9.3.5, 11.3.5]{Du}.
\end{myItem}

\begin{myItem}
Let $ X $ be a Polish space, and $ K $ a metrizable compact
topological space. Denote by $ C(K,X) $ the space of all continuous
maps $ K \to X $ with the metric $ \dist (f,g) = \sup_{a \in K} \dist
( f(a), g(a) ) $; equivalent metrics on $ X $ induce equivalent
metrics on $ C(K,X) $, and $ C(K,X) $ is a Polish space. See
\cite[4.19]{Ke} or \cite[2.4.3]{Sri}.
\end{myItem}

\begin{myItem}
Let $ X $ be a Polish space. Denote by $ \K (X) $ the topological
space of all compact subsets of $ X $, equipped with the
\emph{Vietoris topology.} It is the topology that corresponds to
Hausdorff metric
\[
\dist (K_1,K_2) = \sup_{x\in X} | \dist(x,K_1) - \dist(x,K_2) | \, ;
\]
here $ \dist (x,K) = \inf_{y\in K} \dist (x,y) $. We assume that $
\dist (x,y) \le 1 $ for all $ x,y $ (otherwise one can use $ \min \(
1, \dist(x,y) \) $ instead of $ \dist(x,y) $) and let $ \dist
(x,\emptyset) = 1 $. A choice of a (compatible) metric on $ X $
influences the Hausdorff metric on $ \K(X) $ but not the topology on $
\K(X) $. (See \cite[Sect.~4.F]{Ke}, see also \cite[p.~1124]{Be}.)

For every Polish space $ X $, the topological space $ \K (X) $ is
Polish.
(See \cite[4.25]{Ke}.)
\end{myItem}

\begin{myItem}
\label{Beer}
Let $ X $ be a Polish \emph{metric} space. Denote by $ \bF (X) $ the
topological space of
all closed subsets of $ X $, equipped with the \emph{Wijsman
topology;}\footnote{%
 The topology of $ X $ does not determine uniquely the Wijsman
 topology on $ \bF(X) $; equivalent metrics on $ X $ may lead to
 different topologies on $ \bF(X) $.}
it is the unique metrizable topology such that for every $
F, F_1, F_2, \dots \in \bF (X) $
\[
F_n \to F \quad \text{if and only if} \quad \forall x \in X \; \dist
(x,F_n) \to \dist (x,F)
\]
provided that $ \dist (x,y) \le 1 $ for all $ x,y $ (otherwise we use
$ \min \( 1, \dist(x,y) \) $ instead of $ \dist(x,y) $); still, $
\dist (x,\emptyset) = 1 $. Any sequence $ (x_n) $ dense in $ X $ gives
rise to a compatible\footnote{%
 Not complete in general.}
metric on $ \bF (X) $, say,
\[
\dist (F_1,F_2) = \max_n \frac1n | \dist (x_n,F_1) - \dist (x_n,F_2) |
 \, .
\]
(See \cite[pp.~25--26]{BK} and \cite[Sect.~3]{Be}.)

Theorem (Beer \cite{Be}).
For every Polish metric space $ X $, the topological space $ \bF (X) $
is Polish.
\end{myItem}

\begin{myItem}\label{1.14}
Note that $ \bF (X) $ is partially ordered (by inclusion) and the
corresponding set of pairs, $ \{ (F_1,F_2) : F_1 \subset F_2 \} $, is
closed in $ \bF(X) \times \bF(X) $.
Note also the \emph{sandwich argument:} if $ E_n \subset F_n \subset
G_n $, $ E_n \to F $ and $ G_n \to F $ then $ F_n \to F $. It follows
immediately from monotonicity of the correspondence between a set $ F
\in \bF(X) $ and its distance function $ x \mapsto \dist(x,F) $. The
same holds for $ \K(X) $.
\end{myItem}

\begin{lemma}\label{1.9}
Let $ X $ be a Polish metric space and $ f : X \times X \to X $ a
continuous function. Then the following condition on a closed set $ F
\subset X $ selects a \emph{closed} subset of $ \bF (X) $:
\[
\forall x,y \in F \;\; f(x,y) \in F \, .
\]
\end{lemma}

\begin{proof}
Let $ F_n $ satisfy the condition, $ F_n \to F $, and $ x,y \in F
$. We have $ \dist (x,F_n) \to 0 $, $ \dist (y,F_n) \to 0 $. Take $
x_n, y_n \in F_n $ such that $ x_n \to x $, $ y_n \to y $. Then $ F_n
\ni f (x_n,y_n) \to f (x,y) $, therefore $ \dist ( f(x,y), F_n ) \to 0
$. So, $ f(x,y) \in F $.
\end{proof}

The binary operation $ (x,y) \mapsto f(x,y) $ may be replaced with a
unary operation $ x \mapsto f(x) $, or even a constant ($ 0 $-ary
operation) $ f \in X $. Note that a closed subset $ F $ of a separable
Hilbert space $ H $ is a \emph{subspace} if and only if it is closed
under the binary operation $ (x,y) \mapsto x+y $, every one of unary
operations $ x \mapsto ax $, and contains $ 0 $.\footnote{%
 Well, you may reduce the list.}
The conclusion follows.

\begin{corollary}\label{1.10}
Let $ H $ be a separable Hilbert space. The set $ \bL (H) $ of all
(closed linear)
subspaces of $ H $, equipped with the Wijsman topology, is a Polish
space.
\end{corollary}

\begin{lemma}\label{1.11}
\begin{sloppypar}
The following two conditions are equivalent for all $ L, L_1, L_2,
\dots \in \bL (H) $:
\end{sloppypar}

(a) $ L_n \to L $ in $ \bL (H) $;

(b) $ P_n x \to P x $ in $ H $ for every $ x \in H $; here $ P_n $ is
the orthogonal projection onto $ L_n $, and $ P $ --- onto $ L $.
\end{lemma}

\begin{proof}
Condition (a) means that $ \| P_n x - x \| \to \| Px - x \| $ for all
$ x $.

(b) \imp (a):
$ P_n x \to Px $, therefore $ \| P_n x - x \| \to \| Px - x \| $.

(a) \imp (b):
We have $ \| P_n x - x \| \to \| Px - x \| $ and $ \| P_n P x - P x \|
\to 0 $. Note that $ \| y - x \|^2 = \| y - P_n x \|^2 + \| P_n x - x
\|^2 $ for all $ y \in L_n $; in particular,
\[
\| P_n Px - x \|^2 = \| P_n Px - P_n x \|^2 + \| P_n x - x \|^2 \, .
\]
However, $ P_n P x \to P x $, thus $ \| P_n P x - x \|^2 \to \| Px
- x \|^2 $ and $ \| P_n Px - P_n x \|^2 = \| Px - P_n x \|^2 + o(1)
$. So, $ \| Px - x \|^2 + o(1) = \| Px - P_n x \|^2 + \| Px - x \|^2 +
o(1) $, thus $ \| Px - P_n x \|^2 \to 0 $.
\end{proof}

\begin{corollary}\label{1.18}
Let $ (\Om,\F,P) $ be a separable probability space. The set of all
sub-\sif s of $ \F $, equipped with the Wijsman topology, is a Polish
space. (Each \sif\ must contain all sets of measure zero, and is
treated as a closed subset of $ \MALG (\Om,\F,P) $.\footnote{%
 However, it may be also treated as a subspace of $ L_2 (\Om,\F,P) $,
 see \eqref{2a2}.}%
)
\end{corollary}

The proof of the latter is left to the reader. (Similar to \ref{1.10},
but Boolean operations are used instead of linear operations.)

\subsection{Polish groups, their actions, orbits}

\begin{definition}
(a) A \emph{Polish group} is a topological group whose topological
space is Polish.

(b) Let $ G $ be a Polish group. A \emph{Polish $ G $-space} is a
Polish topological space $ X $ equipped with a continuous map $ (g,x)
\mapsto g \cdot x $ from $ G \times X $ to $ X $ such that $ 1 \cdot x
= x $ and $ (gh) \cdot x = g \cdot ( h \cdot x ) $ for all $ x \in X $
and $ g,h \in G $.

(c) Let $ G $ be a Polish group. A \emph{Polish metric $ G $-space} is
a Polish metric space $ X $ equipped with a map $ (g,x) \mapsto g
\cdot x $ satisfying (b) and in addition, $ \dist ( g \cdot x, g \cdot
y ) = \dist (x,y) $ for all $ x,y \in X $, $ g \in G $.

(d) Let $ X $ be a Polish $ G $-space. Its \emph{orbit equivalence
relation} $ E_G $ is defined by
\[
(x,y) \in E_G \quad \text{if and only if} \quad \exists g \in G \;\; g \cdot x =
y \, .
\]
The \emph{orbit} of $ x $ is its equivalence class $ [x]_G = G \cdot x
$.
\end{definition}

From now on, the phrase `Let $ X $ be a Polish $ G $-space' means
`Let $ G $ be a Polish group and $ X $ a Polish $ G $-space'. The same
for Polish metric $ G $-spaces.

\begin{conjecture}(Topological Vaught Conjecture)
Let $ X $ be a Polish $ G $\nobreakdash-space.
If $ X $ has uncountably many orbits then $ X $ contains an
uncountable \emph{closed} set $ C $ with no equivalent points (that
is, if $ y = g \cdot x $ and $ x,y \in C $ then $ x = y $).
(See \cite[Sect.~6.2 and 3.3]{BK}; \cite[Sect.~5.13]{Sri}.)
\end{conjecture}

The subset $ E_G $ of $ X \times X $ need not be a Borel set (see
\cite[Sect.~3.2]{BK}), even though every orbit is a Borel subset of $
X $ (see \cite[2.3.3]{BK}). However, $ E_G $ is an analytic subset of
$ X \times X $, since it is the projection to $ X \times X $ of the
closed set $ \{ (g,x,y) : g \cdot x = y \} \subset G \times X \times X
$. All analytic sets have the Baire property, and are universally
measurable (that is, measurable w.r.t.\ every Borel measure); see
\cite[4.3.1, 4.3.2]{Sri} or \cite[21.6, 21.10, Sect.~29.B]{Ke}, see
also \cite[13.2.6]{Du}.

\begin{lemma}
Let $ X $ be a Polish $ G $-space. Then there exists a function $ g_1
: E_G \to G $ Baire measurable, universally measurable, and such that
\[
g_1 (x,y) \cdot x = y \qquad \text{for all } (x,y) \in E_G \, .
\]
\end{lemma}

\begin{proof}
Follows from Von Neumann's theorem on measurable selection, see
\cite[5.5.3]{Sri}, \cite[29.9]{Ke}, see also \cite[13.2.7]{Du}.
\end{proof}

\begin{definition}\label{1.12}
Let $ X $ be a Polish $ G $-space.

(a) A \emph{$ G $-invariant set} in $ X $ is $ A \subset X $ such that
$ \forall g \in G \; \forall x \in A \; g \cdot x \in A $.

(b) A \emph{$ G $-invariant function} on $ X $ is a map $ f : X \to Z
$ ($ Z $ being an arbitrary set) such that $ f ( g \cdot x ) = f(x) $
for all $ g \in G $, $ x \in X $.

(c) A \emph{selector}\footnote{%
 Or `section'.}
for $ E_G $ is a $ G $-invariant map $ f : X \to
X $ such that $ (x,f(x)) \in E_G $ for all $ x \in X $.\footnote{%
 Of course, $ G $-invariance of $ f $ means $ f ( g \cdot x ) = f(x)
 $, not $ f ( g \cdot x ) = g \cdot f(x) $.}

(d) A \emph{transversal}\footnote{%
 Or `cross-section'.}
for $ E_G $ is a set $ T \subset X $ such that for every $ x \in X $
 the set $ [x]_G \cap T $ contains exactly one point.

(e) The orbit equivalence relation $ E_G $ is called
\emph{smooth}\footnote{%
 Or `tame'.}
if there exist a Polish space $ Z $ and a ($ G $-invariant) Borel
function $ \theta : X \to Z $ such that for all $ x,y \in X $
\[
(x,y) \in E_G \quad \text{if and only if} \quad \theta(x) = \theta(y) \, .
\]
\end{definition}

For (c), (d) see \cite[12.15]{Ke}, \cite[p.~186]{Sri}. For (e) see
\cite[18.20]{Ke}, \cite[0.1]{Hj98}.

\begin{theorem}\label{1.15}
Let $ X $ be a Polish $ G $-space, then the following four conditions
are equivalent.

(a) $ E_G $ is smooth;

(b) there exists a Borel selector;

(c) there exists a Borel transversal;

(d) there exist $ G $-invariant Borel sets $ A_1, A_2, \dots \subset X
$ that separate orbits in the sense that for all $ x,y \in X $
\[
(x,y) \in E_G \quad \text{if and only if} \quad \forall n \, \( x \in
A_n \equiv y \in A_n \) \, .
\]
\end{theorem}

See \cite[5.6.1: Burgess theorem]{Sri} and \cite[18.20(i,iii)]{Ke}. The
transition (c) \imp (b) uses the Lusin separation theorem for analytic
sets (see \cite[28.1]{Ke} or \cite[4.4.1]{Sri}) in a way similar to
(but simpler than) \cite[4.4.5]{Sri}.

If $ E_G $ is smooth then $ E_G $ evidently is a Borel subset of $ X
\times X $. The converse is wrong; for a counterexample, consider the
natural action of the additive group $ G = \Q $ of rational numbers on
the space $ X = \R $ of real numbers.

A Borel selector $ f $ maps bijectively the set of orbits onto the
Borel set $ T_f = \{ x \in X : f(x) = x \} $. Let $ \theta : X \to Z $
be as in \ref{1.12}(e), then the restriction $ \theta |_{T_f} $ maps
bijectively $ T_f $ onto $ \theta(X) $, and is a Borel map. It follows
that $ \theta(X) $ is a Borel set (see
\cite[15.2]{Ke}). Therefore $ Z $ and $ \theta $ may be chosen
so that $ \theta(X) = Z $ (see \cite[13.1]{Ke}; see also
\cite[2.27(iii)]{Hj}).

\begin{definition}\label{1.153}
A Polish $ G $-space $ X $ will be called \emph{ergodic,}\footnote{%
 Not a standard terminology.}
if every $ G $\nobreakdash-invariant set $ A \subset X $ having the
Baire property is either meager or comeager.
\end{definition}

\begin{theorem}\label{1.154}
The following conditions are equivalent for every Polish $ G $-space $
X $.

(a) $ X $ is ergodic.

(b) Every nonempty open $ G $-invariant set is dense.

(c) There exists a dense orbit (at least one).

(d) The union of all dense orbits is a dense $ G_\de $-set.\footnote{%
 A $ G_\de $ set is, by definition, the intersection of a sequence of
 open sets.}

(e) For every Polish space $ Z $, every $ G $-invariant Baire
measurable function $ f : X \to Z $ is constant on a comeager set.

(f) For every Polish space $ Z $ and every Baire measurable function $
f : X \to Z $, if $ f $ is almost invariant in the sense
that\footnote{%
 Recall Sect.~1a about ``$ \forall^* $''. Of course, saying ``$
 \forall g \, \forall^* x $'' we mean that the comeager set of these $
 x $ may depend on $ g $.}
\[
\forall g \in G \;\; \forall^* x \in X \;\; f(g\cdot x) = f(x) \, ,
\]
then $ f $ is almost constant in the sense that
\[
\exists z \in Z \;\; \forall^* x \in X \;\; f(x) = z \, .
\]
\end{theorem}

See also \cite[3.2, 3.4]{Hj}.

\begin{proof}
(a) \imp (b): a meager set cannot have interior points; a comeager set
must be dense.

(c) \imp (b): the open set must contain the dense orbit.

(d) \imp (c): trivial.

(e) \imp (a): just take $ Z $ of only two points.

(f) \imp (e): trivial.

(b) \imp (a): Given a $ G $-invariant set $ A \subset X $ that
has the Baire property, we need a $ G $-invariant open set $ U \subset
X $ such that $ U \setminus A $ and $ A \setminus U $ are
meager. However, the canonical construction of $ U $ given by
\cite[8.29]{Ke} is evidently $ G $\nobreakdash-invariant.

(b) \imp (d): It is easy to see that the union of all dense orbits is
the intersection of all nonempty open $ G $-invariant sets. Without
affecting the intersection we can restrict ourselves to a sequence of
such sets (namely, $ G $-saturations of sets of a countable base).

(a) \imp (e): For every open $ U \subset Z $ the set $ f^{-1} (U) $ is
either meager or comeager. Consider the union $ V $ of all open sets $
U \subset Z $ such that $ f^{-1} (U) $ is meager. The set $ f^{-1} (V)
$ is meager (since we may restrict ourselves to a countable base). It
remains to prove that $ Z \setminus V $ contains only a single
point. However, for two disjoint open sets $ U_1, U_2 $ the sets $
f^{-1} (U_1) $, $ f^{-1} (U_2) $, being disjoint, cannot be comeager
simultaneously.

(b) \imp (f): Similarly to ``(b) \imp (a) \imp (e)'' we have (b) \imp
($ \text{a}' $) \imp (f), where ($ \text{a}' $) says that every
\emph{almost} invariant set 
$ A \subset X $ (that is, such that $ A $ and $ g\cdot A $ differ by a
meager set) having Baire property is either meager or comeager. The
construction used when proving ``(b) \imp (a)'' still works, giving an
\emph{invariant} $ U $ for an \emph{almost} invariant $ A $.
\end{proof}

\begin{corollary}\label{1.155}
Let $ X $ be an ergodic Polish $ G $-space such that $ E_G $ is
smooth. Then there exists a comeager orbit.
(See also \cite[4.3]{Hj98}, \cite[3.3, 3.5]{Hj}.)
\end{corollary}

\subsection{Constructing Polish groups and Polish $ G $-spaces}

Recall constructions of 1b: $ C(K,X) $, $ \K(X) $, $ \bF(X) $, $ \MALG
(\Om,\F,P) $ and others.

\begin{myItem}
Let $ X $ be a Polish $ G $-space. Then $ \K (X) $ is a Polish $ G
$-space.

See \cite[item (ii) in Sect.~2.4]{BK}. Here and henceforth we do not
specify the action, having in mind the evident, natural action.
\end{myItem}

\begin{myItem}\label{1.28}
Let $ X $ be a Polish $ G $-space, and $ K $ a metrizable compact
topological space. Then $ C(K,X) $ is a Polish $ G $-space.

The proof is left to the reader.
\end{myItem}

Let $ X $ be a Polish \emph{metric} space. Denote by $ \Iso (X) $ the
group of \emph{isometric} invertible transformations $ \al : X \to X
$.\footnote{%
 For some $ X $ the group $ \Iso(X) $ may be small, say, the unit
 only. For some other $ X $, however, $ \Iso(X) $ is quite large.}
Equip $ \Iso (X) $ with the unique metrizable topology such that for
every $ \al, \al_1, \al_2, \dots \in \Iso (X) $
\[
\al_n \to \al \quad \text{if and only if} \quad \forall x \in X \;\;
\al_n (x) \to \al (x) \, .
\]
Any sequence $ (x_n) $ dense in $ X $ gives rise to a compatible
metric on $ \Iso (X) $, say,
\[
\dist (\al,\be) = \max_n \frac1n \dist \( \al(x_n), \be(x_n) \) \, .
\]

\begin{myItem}\label{1.16}
Theorem.
Let $ X $ be a Polish metric space, then

(a) $ \Iso (X) $ is a Polish group,

(b) $ X $ is a Polish metric $ \Iso(X) $-space.

For Item (a) see \cite[item (9) in Sect.~9.B]{Ke}. Item (b) is left to
the reader. (Continuity of $ g \cdot x $ may be checked in $ g $ and $
x $ separately due to \cite[9.14]{Ke}.)
\end{myItem}

\begin{myItem}
\label{1.23}
Let $ X $ be a Polish metric $ G $-space. Then $ \bF(X) $ is a Polish
(topological) $ G $-space.

The proof is left to the reader. (Once again, separate continuity is
enough.)
\end{myItem}

\begin{lemma}\label{1.17}
Let $ X $ be a Polish metric space and $ f : X \times X \to X $ a
continuous map. Then the following condition on $ \al \in \Iso(X) $
selects a \emph{closed} subset of $ \Iso(X) $:
\[
\forall x,y \in X \;\; f \( \al(x), \al(y) \) = \al \( f(x,y) \) \, .
\]
\end{lemma}

\begin{proof}
Similarly to \ref{1.9}, if $ \al_n \to \al $ then $ f \( \al(x), \al(y) \) =
\lim f \( \al_n(x), \al_n(y) \) = \lim \al_n \( f(x,y) \) = \al \( f(x,y) \)
$.
\end{proof}

The same holds for unary operations and constants. Some applications
follow.

\begin{myItem}\label{1.32}
Let $ H $ be a separable Hilbert space. The group $ \U(H) $ of all
isometric linear invertible operators $ H \to H $ is a Polish group;
and $ H $ is a Polish metric $ \U(H) $-space. Also the set $ \bL(H) $
of all (closed linear) subspaces of $ H $ is a Polish $ \U(H)
$-space.
\end{myItem}

\begin{myItem}\label{1.33}
Let $ (\Om,\F,P) $ be a separable probability space. The group $ \Aut
(\Om,\F,P) $ of all measure preserving automorphisms of $ \MALG
(\Om,\F,P) $\footnote{%
 May we say `measure preserving maps $ \Om \to \Om $'? We'll return
 to the question in Sect.~2.}
is a closed subgroup of $ \Iso \( \MALG (\Om,\F,P) \) $, therefore a
Polish group, and $ \MALG (\Om,\F,P) $ is a Polish metric $ \Aut
(\Om,\F,P) $-space. Therefore the set of all sub-\sif s of $ \F $
(recall \ref{1.18}) is a
Polish $ \Aut (\Om,\F,P) $-space. (See also \cite[17.46(i)]{Ke}.)
\end{myItem}

One may also consider transformations that send $ P $ into
equivalent measures (not just into itself). (See
\cite[17.46(ii)]{Ke}.)

\begin{myItem}
Another construction (unrelated to \ref{1.16}, \ref{1.17}).
Let $ X $ be a metrizable \emph{compact} topological space, and $
\Homeo (X) $ the group of all homeomorphisms $ X \to X $. Being
equipped with a natural topology, $ \Homeo (X) $ is a Polish group,
and $ X $ is a Polish $ \Homeo(X) $-space, as well as $ \K(X) $.
(See \cite[item (8) in Sect.~9.B]{Ke} and \cite[4.2]{Hj}.)
\end{myItem}

\section{Preliminaries on filtrations (for specialists in the
classification theory)}
\subsection{Two equivalent languages}

Probability theory speaks two equivalent languages, `pointful' and
`pointless'. In the `pointful' language, a morphism between two
probability spaces $ (\Om_1,\F_1,P_1) $ and $ (\Om_2,\F_2,P_2) $ is
defined as a measure preserving map $ \al : \Om_1 \to \Om_2 $ (or
rather, equivalence class of such maps). Evidently, it generates maps
\begin{equation}\label{2a}
\begin{gathered}
L_2 (\Om_2,\F_2,P_2) \to L_2 (\Om_1,\F_1,P_1) \, , \\
L_0 (\Om_2,\F_2,P_2) \to L_0 (\Om_1,\F_1,P_1) \, , \\
\MALG (\Om_2,\F_2,P_2) \to \MALG (\Om_1,\F_1,P_1) \, .
\end{gathered}
\end{equation}
These maps are continuous, linear (for $ L_2 $ and $ L_0 $), preserve
structural operations (`$ \max $' and `$ \min $' for $ L_2 $ and $ L_0 $;
`$ \cup $' and `$ \cap $' for $ \MALG $), distributions of random
variables (for $ L_2 $ and $ L_0 $) or probabilities of events (for $
\MALG $). Given a map (between $ L_2 $, or $ L_0 $, or $ \MALG $)
having such properties, one can reconstruct the corresponding $ \al $
(it is unique $ \bmod \, 0 $ and exists) provided, however, that the
probability spaces are good enough. Namely, each probability space is
assumed to be a Lebesgue-Rokhlin space, that is, isomorphic ($ \bmod 0
$) to an interval with Lebesgue measure, or a (finite or countable)
set of atoms, or a combination of the two.

The `pointless' language defines a morphism as a map between spaces $
L_2 $ (or $ L_0 $, or $ \MALG $) satisfying the properties mentioned
above. Thus, it avoids any restrictions on probability spaces, except
for separability. Moreover, one can escape $ \Om $ at all, since $
\MALG $ (or $ L_0 $) can be axiomatized.
(See \cite[17.F]{Ke} and \cite[Sect. 1]{BEKSY}; see also \ref{3.2}.)

Both languages describe the same category. An \emph{object} is either
a Lebesgue-Rokhlin probability space $ (\Om,\F,P) $ or, equivalently,
a separable measure algebra $ \MALG $ (as defined in
\cite[17.44]{Ke}); also a linear lattice $ L_2 $ or $ L_0 $ may be
used. A
\emph{morphism} is either $ \Om_1 \to \Om_2 $ or $ \MALG_2 \to \MALG_1
$ (see above). In the `pointful' language, a morphism from the first
object to the second is surely $ \Om_1 \to \Om_2 $. In the
`pointless' language we have no consensus about the direction; some
authors define a morphism from the first object to the second as $
\MALG_1 \to \MALG_2 $. However, I prefer to define the direction
according to $ \Om_1 \to \Om_2 $ in all cases, even if $ \Om_1, \Om_2
$ are implicit. That is, for me \eqref{2a} is a morphism \emph{from
the first object to the second.}

A morphism $ \Om_1 \to \Om_2 $ selects a closed (linear) sublattice in
$ L_2 (\Om_1) $, or $ L_0 (\Om_1) $, or $ \MALG (\Om_1) $, that is,
a topologically closed set, closed under linear (for $ L_2, L_0 $) and
structural operations. Every subset closed in that sense corresponds
to some morphism. Such a subset is a \emph{sub-$\si$-field,} --- just
by definition, if we are speaking the `pointless' language. The
corresponding notion in the `pointful' language is a sub-\sif\ of $
\F $ containing all sets of measure $ 0 $. Or equivalently, it is a
measurable partition (see \cite{Ro}). Every sub-\sif\ gives rise to a
quotient space. The latter is again a Lebesgue-Rokhlin space. If $
(\Om,\F,P) $ is an object and $ \E \subset \F $ a sub-\sif\ then $
(\Om,\E,P|_\E) $ is another object. I treat the latter as a quotient
of the former, since the latter is naturally isomorphic to the
quotient space $ \Om / \E $ equipped with its natural \sif\ and
measure; in that sense, $ (\Om,\E,P|_\E) = (\Om,\F,P) / \E $. Note the
natural morphism from the object to its quotient-object.\footnote{%
 The direction of the morphism shows that the smaller object should
 not be called a sub-object of the larger object (unless you prefer
 the other direction of morphisms).}

\textsc{From now on, every probability space is (by assumption or
construction) a Lebesgue-Rokhlin space.}

Given a probability space $ (\Om,\F,P) $, the set $ \SF (\Om,\F,P) $
of all sub-\sif s $ \E \subset \F $ may be treated as a closed subset
of $ \bL ( L_2 (\Om,\F,P) )  \subset \bF ( L_2 (\Om,\F,P) ) $, or $
\bF ( L_0 (\Om,\F,P) ) $, or $ \bF
( \MALG (\Om,\F,P) ) $. Probably it is all the same, but anyway, I
prefer the first option,
\begin{equation}\label{2a2}
\SF (\Om,\F,P) \subset \bL ( L_2 (\Om,\F,P) ) \, .
\end{equation}
Thus, $ \SF (\Om,\F,P) $ is a Polish space. Every
monotone sequence $ (\F_n ) $ of \sif s converges to some \sif\ $ \F $;
if $ (\F_n) $ decreases then $ \F $ is the intersection; if $ (\F_n)
$ increases then $ \F $ is the least \sif\ containing the union.
The least among all sub-\sif s is the trivial \sif\ $ \Zero_\SF $
containing sets of measure $ 0 $ or $ 1 $ only; the greatest one is $
\One_\SF = \F $.

Let $ \E $ be a sub-\sif, then $ \E $-measurable random variables are
a subspace $ L_0 (\E) $ of the space $ L_0 = L_0 (\F) $. Also, $ L_2
(\E) = L_0 (\E) \cap L_2 $. The orthogonal projection $ L_2 \to L_2
(\E) $ has a special notation
\[
f \mapsto \cE f \E \, .
\]
It has a continuous extension (evidently unique) to an operator $ L_1
\to L_1 (\E) $, still denoted by $ f \mapsto \cE f \E $. The random
variable $ \cE f \E $ is called the \emph{conditional expectation} of
$ f $ w.r.t.\ $ \E $. Note that
\begin{multline}\label{2a3}
\E_n \to \E \text{ in $ \SF(\Om,\F,P) $} \qquad \text{if and only if}
 \\
\forall f \in L_2 (\Om,\F,P) \;\; \cE f {\E_n} \to \cE f \E \text{ in
  $ L_2 (\Om,\F,P) $}
\end{multline}
by Lemma \ref{1.11}. An equivalent condition:
\[
\forall A \in \F \quad \cP A {\E_n} \to \cP A \E \text{ in
  $ L_2 (\Om,\F,P) $} \, ;
\]
the conditional probability $ \cP A \E $ is, by definition, $ \cE{
\One_A }{ \E } $, where $ \One_A $ is the indicator, $ \One_A (\om) =
1 $ for $ \om \in A $, otherwise $ 0 $. See also \cite[10.2.9]{Du}.

The conditional distribution $ P_{f|\E} $ is an element of $ L_0
(\Om,\F,P; \Prob(\R) ) $\footnote{%
 Recall \ref{1.9a}, \ref{1.10a}.}
such that $ P_{f|\E} (\cdot) (A) = \cP A \E $
for every Borel set $ A \subset \R $, or equivalently, $ \int_{\R} \phi \,
dP_{f|\E} (\cdot) = \cE{ \phi\circ f }{ \E } $ for every bounded Borel
function $ \phi : \R \to \R $. Such $ P_{f|\E} $ exists, is unique,
and is $ \E $-measurable.\footnote{%
 You may derive it from \ref{3.2}. See also \cite[10.2.2]{Du}.}
The same for $ f \in L_0 (\Om,\F,P; X) $ and $ P_{f|\E} \in L_0
(\Om,\F,P; \Prob(X)) $,
where $ X $ is a Polish space.

\subsection{Filtrations}

\begin{definition}\label{2.3}
(a) A \emph{filtration}\footnote{%
 The measure $ P $ is
 given. Sometimes a filtration is considered rather on a \emph{measure
 type space,} that is, only an equivalence class of the measure is
 given. Maybe I should say `filtered probability space' instead.}
(on a probability space $ (\Om,\F,P) $) is a family $
(\F_t)_{t\in[0,\infty)} $ of sub-\sif s $ \F_t \subset \F $ satisfying
 $ \F_s \subset \F_t $ whenever $ 0 \le s \le t < \infty $.

(b) A filtration is called \emph{continuous} if the function $ t
\mapsto \F_t $ is continuous on $ [0,\infty) $. The same for
`right-continuous' and `left-continuous'.
\end{definition}

One-sided limits $ \F_{t-} $, $ \F_{t+} $ arise naturally. Note that $
\F_{0+} $ need not be the trivial \sif\ $ \Zero_\SF $, and $ \F_\infty
$ need not be the whole $ \One_\SF = \F $.

\begin{definition}\label{2.4}
(a) A \emph{random process} (on a probability space $ (\Om,\F,P) $) is
a family $ (f_t)_{t\in[0,\infty)} $ of random variables $ f_t \in L_0
$.\footnote{%
 Some non-equivalent definitions are in use; say, one may stipulate a
measurable function $ [0,\infty) \times \Om \to \R $.}

(b) The \emph{natural filtration} $ (\F_t) $ of a random process $
(f_t) $ (or the \emph{filtration generated by} a random process $ (f_t) $) is
defined as follows: for each $ t $, $ \F_t $ is the sub-\sif\
generated by $ (f_s)_{s\in[0,t]} $, that is, the least sub-\sif\ such
that $ f_s \in L_0 (\F_t) $ for all $ s \in [0,t] $.

(c) A random process $ (f_t) $ is a \emph{martingale} w.r.t.\ a filtration
$ (\F_t) $, if
\[
f_t \in L_1 (\F_t) \quad \text{and} \quad f_s = \cE{ f_t }{ \F_s }
\]
whenever $ 0 \le s \le t < \infty $.

(d) A filtration $ (\E_t) $ is \emph{immersed} into another filtration
$ (\F_t) $, if every martingale w.r.t.\ $ (\E_t) $ is a martingale
w.r.t.\ $ (\F_t) $.
\end{definition}

Condition (d) may be written as $ \cE{ \cE{ f }{ \E_t } }{ \F_s } = \cE{
f }{\E_s } $ whenever $ 0 \le s \le t < \infty $, or in operator form,
$ \F_s \E_t = \E_s $ if we identify sub-\sif s with their operators of
conditional expectation. An implication of (d): $ \E_t \F_s = \F_s
\E_t $ (since $ \E_t \F_s = (\F^*_s \E^*_t)^* = (\F_s \E_t)^* = \E^*_s
= \E_s = \F_s \E_t $), thus, all the operators belong to a commutative
algebra. Another implication: $ \E_s = \F_s \cap \E_t $, especially, $
\E_t = \F_t \cap \E_\infty $. An immersed filtration $ (\E_t) $ is
uniquely determined by $ \E_\infty $ (and $ (\F_t) $).\footnote{%
 Another equivalent form of (d): $ \E_\infty $ and $ \F_t $ are
 conditionally independent, given $ \E_t $.}

An isomorphism between two filtrations is defined in such a way that
interrelations between all $ \F_t $ are relevant, but interrelations
between $ \F_t $ and $ \F $ are not.

\begin{definition}\label{2.5}
Let $ (\F_t) $ be a filtration on a probability space $ (\Om,\F,P) $
and $ (\F'_t) $ a filtration on another probability space $
(\Om',\F',P') $.

(a) An \emph{isomorphism} between the two filtrations is
an \emph{invertible} morphism between quotient spaces $
(\Om,\F,P)/\F_\infty $ and $ (\Om',\F',P')/\F'_\infty $ such that the
image of $ \F_t $ is $ \F'_t $ for every $ t $.

(b) A \emph{morphism} from $ (\F_t) $ to $ (\F'_t) $ is an isomorphism
between $ (\F'_t) $ and a filtration immersed into $ (\F_t) $.
\end{definition}

See also \cite[Sect.~1,2]{BsE}, \cite[Sect.~1]{Wat}.

\subsection{Few examples of random processes and filtrations}

\begin{myItem}\label{2.7}
Sub-\sif s $ \E_1, \E_2, \dots $ are called \emph{independent,} if $ P
( A_1 \cap A_2 \cap \dots ) = P(A_1) P(A_2) \dots $ for all $ A_1 \in
\E_1, A_2 \in \E_2, \dots \, $
Random variables $ f_1, f_2, \dots $ are called independent, if
they generate independent sub-\sif s.
\end{myItem}

\subsubsection{Poisson process}

Take a sequence of independent random variables\footnote{%
 No need to specify the probability space and the choice of measurable
 functions $ \xi_k $ on it, since the filtration $ (\F_t) $
 constructed below is determined uniquely up to isomorphism.}
$ \xi_1, \xi_2, \dots $, each distributed exponentially, namely, for
every $ k $
\[
P ( \xi_k \le c ) = 1 - e^{-c} \quad \text{for all } c \in [0,\infty)
\, .
\]
Define a random process $ (f_t) $ by
\[
f_t(\om) = \max \{ k : \xi_1(\om) + \dots + \xi_k(\om) \le t \} \, ,
\]
then increments $ f_{t_n}-f_{t_{n-1}}, \dots, f_{t_1}-f_{t_0} $ are
independent whenever $ 0 \le t_0 \le \dots \le t_n < \infty $, and
have Poisson distributions:
\[
\Pr{ f_t - f_s = k } = \frac1{k!} (t-s)^k e^{-(t-s)} \quad \text{for }
k = 0,1,2,\dots \text{ and } s \le t \, .
\]
Such $ (f_t) $ is called Poisson process. (See also \cite[Problem
12.1.12]{Du}.) Its natural filtration $
(\F_t) $ is continuous. Every $ \F_t $ has a single atom (of
probability $ e^{-t} $) and a nonatomic part. Here are two examples
of martingales w.r.t.\ $ (\F_t) $:
\[
M_t = f_t - t \, ; \quad N_t = M_t^2 - t \, .
\]
Another process $ (g_t) $ defined by
\[
g_t = \entier \( \tfrac12 f_t ) = \max \{ k : \xi_1 + \dots + \xi_{2k}
\le t \}
\]
has its natural filtration immersed into $ (\F_t) $. On the other
hand, the `slowed down' filtration $ (\E_t) $ defined by
\[
\E_t = \F_{t/2} \quad \text{for } t \in [0,\infty)
\]
is \emph{not} immersed into $ (\F_t) $, even though $ \E_t \subset
\F_t $ for all $ t $. In fact, the three filtrations are pairwise
non-isomorphic.

\subsubsection{Brownian motion}

Take a sequence of independent random variables $ \xi_1, \xi_2, \dots
$, each distributed normally, namely, for every $ k $
\[
P ( \xi_k \le c ) = (2\pi)^{-1/2} \int_{-\infty}^c e^{-u^2/2} \, du
\qquad \text{for all } c \in \R \, . 
\]
The (closed linear) subspace $ G \subset L_2 (\Om,\F,P) $ spanned by $
(\xi_k) $ consists of normally distributed random variables, and
independence is equivalent to orthogonality within $ G $ (see
\cite[9.5.14]{Du}). Take any linear isometric operator
\[
U : L_2 (0,\infty) \to G
\]
and define a random process $ (B_t) $ by
\[
B_t = U \( \One_{(0,t)} \) \quad \text{for } t \in [0,\infty) \, ;
\]
here $ \One_{(0,t)} \in L_2 (0,\infty) $ is the indicator of the
interval $ (0,t) $.
Increments $ X_{t_n}-X_{t_{n-1}}, \dots, X_{t_1}-X_{t_0} $ are
independent whenever $ 0 \le t_0 \le \dots \le t_n < \infty $, and
have normal distributions:
\[
\Pr{ B_t - B_s \le \sqrt{t-s} \, c } = (2\pi)^{-1/2} \int_{-\infty}^c
e^{-u^2/2} \, du \quad \text{for all } c \in \R, s \le t \, .
\]
Such $ (B_t) $ is called Brownian motion. (See also
\cite[12.1.5]{Du}.) Its natural filtration $
(\F_t) $ is continuous. Every $ \F_t $ is non-atomic (except for $
\F_0 $). Here are two examples of martingales w.r.t.\ $ (\F_t) $:
\[
M_t = B_t \, ; \quad N_t = M_t^2 - t \, .
\]
The `slowed down' filtration $ (\F_{t/2})_{t\in[0,\infty)} $ is
isomorphic to $ (\F_t)_{t\in[0,\infty)} $, but not immersed.

\section{Isomorphism of probability spaces\\
 as an orbit equivalence relation}
\subsection{Rokhlin's theory revisited}

Rokhlin's theory \cite{Ro} gives us, first, a classification of
probability spaces\footnote{%
 Do not forget our convention proclaimed in Sect.~2a: every
 probability space is (by assumption or construction) a
 Lebesgue-Rokhlin space.}
(up to isomorphisms $ \bmod \, 0 $), and second, a classification of
sub-\sif s, or rather, pairs $ \( (\Om,\F,P), \E \) $, where $
(\Om,\F,P) $ is a probability space and $ \E \subset \F $ a
sub-\sif.\footnote{%
 Containing all $ P $-negligible sets (recall Sect.~2a).}
In other words, the second is a classification of morphisms\footnote{%
 That is, measure preserving maps (recall Sect.~2a).}
$ \al : (\Om_1,\F_1,P_1) \to (\Om_2,\F_2,P_2) $, or rather, triples $ \(
(\Om_1,\F_1,P_1),  (\Om_2,\F_2,P_2), \al \) $. These results are
formulated below (without proofs). Orbit equivalence relations do not
appear in Sect.~3a. A part of Rokhlin's theory (roughly, the
classification of probability spaces, but not morphisms) is presented
also in \cite{Ha}, \cite{Ru}.

Let us choose once and for all our `favorites': an uncountable
standard Borel space $
(\Omfav,\Ffav) $, say, $ (0,1) $ with the Borel \sif; and a nonatomic
probability measure $ \Pfav $ on $ (\Omfav,\Ffav) $, say, the Lebesgue
measure on $ (0,1) $; and a sequence of different points $ \om_1,
\om_2, \dots \in \Om_0 $, say, $ \om_k = 2^{-k} $. We introduce the
set $ \Mfav $, consisting of all probability measures $ \mu $ on $
(\Omfav,\Ffav) $ of the form
\[
\mu (A) = m_0 \Pfav (A) + \sum_{k=1}^\infty m_k \One_A (\om_k) \, ,
\]
where $ m_k \in [0,1] $, $ m_0 + m_1 + m_2 + \dots = 1 $ and in
addition, $ m_1 \ge m_2 \ge \dots $ In other words, $ \mu $ is a
convex combination of atoms at $ \om_k $ (with decreasing weights $
m_k $), and the nonatomic measure $ \Pfav $.

\begin{theorem}\label{3.1}
Every probability space is isomorphic ($ \bmod \, 0 $) to $
(\Omfav,\mu) $ for one and only one measure $ \mu \in \Mfav $.
\end{theorem}

See \cite[Sect.~2.4]{Ro}. Of course, the \sif\ of $ (\Omfav,\mu) $ is
$ \Ffav $ completed by adding $ \mu $-negligible sets. We'll usually
suppress \sif s in the notation like $ (\Om,P) $ when it is just the
natural \sif\ on $ \Om $, completed by $ P $-negligible sets.

\begin{note}
All nonatomic probability spaces are isomorphic to $ (\Omfav,\Pfav)
$.
\end{note}

Probability spaces are thus classified. In order to classify morphisms
(in other words, sub-\sif s, or measurable partitions) we introduce
another set $ \M_2 $, consisting of all probability measures $ \mu $
on $ (\Omfav,\Ffav) \times (\Omfav,\Ffav) $ of the form
\[
\mu ( A \times B ) = \int_A \mu_a (B) \, d\mu_1(a) \, ,
\]
where $ \mu_1 \in \Mfav $, and $ \mu_a \in \Mfav $ for $ \mu_1
$-almost all $ a $, and the map $ a \mapsto \mu_a $ is $ \mu_1
$-measurable. (In other words: both the marginal distribution $ \mu_1
$ and conditional distributions $ \mu_a $ must belong to the class $
\Mfav $.)

Especially, if $ \mu_a $ is nonatomic for all $ a $, then $ \mu_a =
\Pfav $ and so, $ \mu = \mu_1 \otimes \Pfav $.

Every measure $ \mu \in \M_2 $ determines a morphism
\[
( \Omfav \times \Omfav, \mu ) \xrightarrow{\mathrm{pr}} (\Omfav,\mu_1) \,
,
\]
just the projection, $ (a,b) \mapsto a $.

\begin{theorem}\label{3.2}
For every morphism $ (\Om_1,P_1) \xrightarrow{\al} (\Om_2,P_2) $ of
probability spaces there exist $ \mu \in \M_2 $ and isomorphisms $
\be, \ga $ such that the diagram
\[
\xymatrix{
(\Om_1,P_1) \ar[d]_\al \ar@{<->}[r]^-{\be} & (\Omfav \times \Omfav, \mu)
 \ar[d]^{\mathrm{pr}} \\
(\Om_2,P_2) \ar@{<->}[r]^-{\ga} & (\Omfav,\mu_1)
}
\]
is commutative.
\end{theorem}

See \cite[Sect.~4.1]{Ro}. Unlike \ref{3.1}, $ \mu $ is (in general)
not uniquely determined by $ \al $.\footnote{%
 I mean, by the triple $ \( (\Om_1,\F_1,P_1), (\Om_2,\F_2,P_2), \al \)
 $.}
If $ \mu_1 $ is purely atomic and all atoms have different
probabilities, then $ \mu $ is unique. However, if there are atoms of
equal probability, say, $ \mu_1 ( \{\om_1\} ) = \mu_1 ( \{\om_2\} ) $
(that is, $ m_1 = m_2 $ for $ \mu_1 $), then conditional measures $
\mu_a $ for $ a = \om_1 $ and $ a = \om_2 $ may be interchanged
(correcting isomorphisms $ \be, \ga $ accordingly). Though, uniqueness
may persist, if these two conditional measures are equal. Similarly,
if $ \mu_1 $ has a nonatomic part, corresponding conditional measures
$ \mu_a $ can be interchanged (but may happen to be equal). This is
why $ \mu \in \M_2 $ cannot be used as an invariant of $ \al $.

Especially, the case of $ \mu = \mu_1 \otimes \Pfav $ is the so-called
\emph{conditionally nonatomic} case. Here, of course, $ \mu $ is
uniquely determined by $ \al $.

We'll return to morphisms after classifying measurable functions.

A measurable function $ f : \Om \to \R $ on a probability space $
(\Om,\F,P) $ determines a morphism $ (\Om,P) \to (\R,P_f) $ where $
P_f $ is the distribution of $ f $ (that is, $ P_f(A) = P \( f^{-1}(A)
\) $). However, an isomorphism between two such functions, $ f_1 :
(\Om_1,P_1) \to \R $ and $ f_2 : (\Om_2,P_2) \to \R $, is defined by a
commutative diagram
\[
\xymatrix{
(\Om_1,P_1) \ar[dr]_{f_1} \ar@{<->}[rr]^-{\al} && (\Om_2, P_2)
 \ar[dl]^{f_2} \\
& (\R,P_f)
}
\]
Only the domain is transformed; the range ($ \R $) is not. The
equality $ P_{f_1} = P_{f_2} (=P_f) $ is necessary but not sufficient
for existence of an isomorphism. Conditional distributions are indexed
by values $ x $ of $ f $ and cannot be interchanged.

We consider the class $ \M_2(P_f) $ of all measures $ \mu $ on $ \R
\times \Omfav $ of the form
\[
\mu (A \times B) = \int_A \mu_x (B) \, d P_f (x)
\]
where $ \mu_x \in \Mfav $ for $ P_f $-almost all $ x \in \R $, and the
map $ x \mapsto \mu_x $ is $ P_f
$\nobreakdash-\hspace{0pt}measurable.\footnote{%
 Thus, $ \M_2 $ introduced before is the union of sets $ \M_2(\mu_1) $
 over all $ \mu_1 \in \Mfav $.}

\begin{theorem}\label{3a3}
For every probability space $ (\Om,\F,P) $ and every measurable
function $ f : \Om \to \R $ there exists one and only one measure $
\mu \in \M_2 (P_f) $ such that there exists an isomorphism $ \al $
such that the diagram
\[
\xymatrix{
(\Om,P) \ar[dr]_{f} \ar@{<->}[rr]^-{\al} && (\R \times \Omfav, \mu)
 \ar[dl]^{\mathrm{pr}} \\
& (\R,P_f)
}
\]
is commutative.
\end{theorem}

The proof is left to the reader (hint: (re)read
\cite[Sect.~4]{Ro}). Note that $ \mu $ is uniquely determined by $ f
$, but $ \al $ is not. Say, if some $ \mu_x $ has two atoms of equal
probability, their interchange influences $ \al $ (but not $ \mu
$). In some sense, $ \mu_x $ describes the multiplicity of the value $
x $ of $ f $.

Note also the conditionally nonatomic case: $ \mu = P_f \otimes \Pfav
$; here, $ \al $ is highly non-unique.

There is nothing special in $ \R $ as the range of $ f $; any other
Polish (or standard Borel) space may be used instead. Especially, $
\Mfav $-valued measurable functions will be used in Theorem \ref{3a4}.

Now we return to morphisms. Recall that the measure $ \mu \in \M_2 $
is not uniquely determined by a morphism $ \al : (\Om_1,P_1) \to
(\Om_2,P_2) $. Note however that the projection $ \mu_1 $ of $ \mu $
is uniquely determined by the morphism, since $ (\Om_2,P_2) $ is
isomorphic to $ (\Omfav,\mu_1) $.

\begin{theorem}\label{3a4}
Let two measures $ \mu', \mu'' $ on $ \Omfav \times \Omfav $ belong to $
\M_2 $, both having the same projection $ \mu_1 \in \Mfav $. Then
existence of isomorphisms $ \al, \be $ making commutative the diagram
\[
\xymatrix{
(\Omfav \times \Omfav,\mu') \ar[d]_{\mathrm{pr}} \ar@{<->}[r]^-{\al} &
 (\Omfav \times \Omfav, \mu'')  \ar[d]^{\mathrm{pr}} \\
(\Omfav,\mu_1) \ar@{<->}[r]^-{\be} & (\Omfav,\mu_1)
}
\]
is equivalent to existence of an isomorphism $ \gamma $ making
commutative the diagram
\[
\xymatrix{
(\Omfav,\mu_1) \ar[dr]_{f'} \ar@{<->}[rr]^-{\gamma} && (\Omfav, \mu_1)
 \ar[dl]^{f''} \\
& \Mfav
}
\]
where a measurable function $ f' : \Omfav \to \Mfav $ is defined by $
f'(a) = \mu'_a $, $ \mu'_a $ being the conditional measure of $ \mu'
$;\footnote{%
 That is, $ \mu'(A\times B) = \int_A \mu'_a (B) \, d\mu_1(a) $.}
the same for $ f'' $ (and $ \mu''_a, \mu'' $).
\end{theorem}

See \cite[Sect.~4.1]{Ro}. Existence of $ \gamma $ may be checked via
the invariant given in Theorem \ref{3a3}. Combining \ref{3a4} and
\ref{3a3} we see that a complete invariant of a morphism $ \al :
(\Om_1,P_1) \to (\Om_2,P_2) $ of probability spaces consists of
\begin{myitemize}
\item
a measure $ \mu_1 \in \Mfav $ on $ \Omfav $, describing the type of the
probability space $ (\Om_2,P_2) $;
\item
a measure $ \mu_2 $ on the space $ \Mfav $ (sic! not a measure
belonging to $ \Mfav $) describing the distribution of the type of the
conditional measure;
\item
a measure $ \mu_3 \in \M_2(\mu_2) $ describing the type of the
conditional measure on $ \Om_1 $, treated as a measurable function on
$ \Om_2 $.\footnote{%
 Thus, $ (\mu_3)_x $ describes the multiplicity of the type $ x \in
 \M^\fav $ of the conditional measure.}
\end{myitemize}

Note the conditionally nonatomic case; here, $ \mu_2 $ is concentrated
at the point $ \Pfav $ of the space $ \Mfav $, and $ \mu_3 = \mu_2
\otimes \mu_1 $.

\subsection{The construction}

Theorem \ref{3.1} shows that probability spaces can be classified by
points of $ \Mfav $, the latter being a Polish space. It should mean
smoothness, as defined by \ref{1.12}(e). To this end we want to embed
Rokhlin's theory into the framework of orbit equivalence relations on
Polish G-spaces. We need to construct a Polish space $ X $ whose
points represent probability spaces, such that all isomorphic classes
are available in $ X $, and all their isomorphisms are available in a
single Polish group $ G $ acting on $ X $.

\begin{sloppypar}
Recall our `favorite' probability space $ (\Omfav,\Pfav)
$.\footnote{%
 In fact, $ (0,1) $ with Lebesgue measure, see Sect.~3a.}
We consider the Polish space
\[
\SF^\fav = \SF (\Omfav,\Pfav)
\]
of all sub-\sif s (see \eqref{2a2}), and the Polish group
\[
\Gfav = \Aut (\Omfav,\Pfav)
\]
of all measure preserving automorphisms (see \ref{1.33}). Note
that
\[
\SF^\fav \text{ is a Polish $ \Gfav $-space,}
\]
since $ \Gfav $ is a closed subgroup of $ \U \( L_2 (\Omfav,\Pfav) \)
$, and $ \SF^\fav \subset \bL \( L_2 (\Omfav,\Pfav) \) $; recall
\ref{1.32}.
\end{sloppypar}

If $ \E_1, \E_2 \in \SF^\fav $ belong to the same orbit then the
corresponding quotient spaces are isomorphic; however, the converse is
wrong. Say, it may happen that $ \E_1, \E_2 \in \SF^\fav $ are both
nonatomic, but $ \E_1 = \Ffav $, $ \E_2 \ne \Ffav $. Then the two
quotient spaces are isomorphic, however, $ \E_1, \E_2 $ belong to
different orbits. For that reason we introduce the set
\begin{equation}\label{3.6a}
X_1 \subset \SF^\fav
\end{equation}
consisting of all \emph{conditionally nonatomic} sub-\sif
s.\footnote{%
 A sub-\sif\ $ \E \in X_1 $ may contain atoms; that is, the
 corresponding measurable partition may contain parts of positive
 probability. However, the conditional measure on (almost) every part
 must be nonatomic. See Sect.~3a. In fact, $ \E $ is conditionally
 nonatomic if and only if there exists a nonatomic sub-\sif\ $ \F $
 independent of $ \E $. It follows easily from \ref{3.2}. See also
 \ref{4.2}.}
Note that for all $ \E_1,\E_2 \in \SF^\fav $
\begin{equation}\label{3.6b}
\text{if } \E_1 \subset \E_2 \text{ and } \E_2 \in X_1 \text{ then }
 \E_1 \in X_1 \, .
\end{equation}

\begin{lemma}\label{3.5}
(a) Every probability space is isomorphic to the quotient space $
(\Omfav,\Pfav) / \E $ for some $ \E \in X_1 $.

(b) For every $ \E_1, \E_2 \in X_1 $ the following two conditions
are equivalent:
\begin{myitemize}
\item quotient spaces $ (\Omfav,\Pfav) / \E_1 $ and $
 (\Omfav,\Pfav) / \E_2 $ are isomorphic;
\item $ \E_2 = g \cdot \E_1 $ for some $ g \in \Gfav $.
\end{myitemize}
\end{lemma}

\begin{proof}
(a) By Theorem \ref{3.1} we may restrict ourselves to probability
spaces $ (\Omfav,\mu) $, where $ \mu \in \Mfav $ is parametrized by $
m_0, m_1, m_2, \dots\, $ We choose disjoint sets $ A_k \subset \Omfav
$ such that $ \Pfav (A_k) = m_k $ and consider the sub-\sif\ $ \E $
generated by all $ A_k $ and all measurable subsets of $ A_0
$. Clearly, $ (\Omfav,\Pfav) / \E $ is isomorphic to $ (\Omfav,\mu)
$. Though, $ \E \notin X_1 $ (unless $ m_0 = 0 $). However, we may
use another probability space $ (\Om_2,P_2) = (\Omfav,\Pfav) \times
(\Omfav,\Pfav) $ and another sub-\sif\ $ \E_2 = \E \otimes \Zero_{\SF}
$; here $ \Zero_{\SF} $ is the trivial \sif\ (consisting of sets of
probability $ 0 $ or $ 1 $ only). Indeed, $ \E_2 $ is conditionally
nonatomic, and $ (\Om_2,P_2) $ is
isomorphic to $ (\Omfav,\Pfav) $, and $ (\Om_2,P_2) / \E_2 $ is
isomorphic to $ (\Omfav,\Pfav) / \E $, therefore, to $ (\Omfav,\mu)
$. 

(b) If $ \E_2 = g \cdot \E_1 $ then $ g $ evidently gives us an
isomorphism between the quotient spaces. On the other hand, let the
quotient spaces be isomorphic, then Theorem \ref{3.2} (specialized for
the conditionally nonatomic case) gives us a commutative diagram
\[
\xymatrix{
(\Omfav,\Pfav) \ar[d] \ar@{<->}[r]^-{\al_1} & (\Omfav \times
 \Omfav,\mu_1 \otimes \Pfav) \ar[d] \ar@{<->}[r]^-{\al_2} &
 (\Omfav,\Pfav) \ar[d] \\
(\Omfav,\Pfav) / \E_1 \ar@{<->}[r] & (\Omfav,\mu_1) \ar@{<->}[r] &
 (\Omfav,\Pfav) / \E_2
}
\]
Combining $ \al_1 $ and $ \al_2 $ we get $ g \in \Gfav $ such that $
\E_2 = g \cdot \E_1 $.
\end{proof}

\begin{note}\label{3.6}
In addition to Item (b) of Lemma \ref{3.5}:

Every isomorphism between $ (\Omfav,\Pfav) / \E_1 $ and $
(\Omfav,\Pfav) / \E_2 $ is induced by some (at least one) $ g \in
\Gfav $.
\end{note}

(A proof is implicitly contained in the proof of Lemma \ref{3.5}.)

\begin{lemma}\label{3.7}
$ X_1 $ is a $ G_\de $-set in $ \SF^\fav $.
\end{lemma}

\begin{proof}
We start with rather general claims; they hold for any probability
space $ (\Om,\F,P) $ and any Polish space $ X $.

\smallskip

\textsc{a.~claim.}
Let $ F \subset X $ be a closed set. Consider such a
function $ \Prob(X) \to \R $:
\[
\Prob(X) \ni \mu \mapsto \mu(F) \in \R \, .
\]
The function is upper semicontinuous.\footnote{%
 That is, $ \{ \mu : \mu(F) < x \} $ is an open subset of $ \Prob(X) $
 for every $ x \in \R $.}

(See \cite[17.20(iii)]{Ke} or \cite[11.1.1(c)]{Du}; see also
\cite[17.29]{Ke}.)

\smallskip

\textsc{b. claim.}
For any $ f \in L_0 (\Om,\F,P; X) $, the map (recall the end of 2a)
\[
\SF (\Om,\F,P) \ni \E \mapsto P_{f|\E} \in L_0 (\Om,\F,P; \Prob(X))
\]
is continuous.

\textit{Proof of Claim.}
Let $ \E,\E_1,\E_2,\dots \in \SF $, $ \E_n \to \E $. It suffices to
find a subsequence $ \E_{n_k} $ such that $ P_{f|\E_{n_k}} \to
P_{f|\E} $ almost everywhere. For any given $ \phi $ (bounded,
continuous, $ X \to \R $) we have $ \cE{ \phi \circ f }{ \E_n } \to
\cE{ \phi \circ f }{ \E } $ in $ L_2 $ by \eqref{2a3}; convergence
almost everywhere holds for some subsequence, however, the subsequence
may depend on $ \phi $. Doing so for a countable set of functions $
\phi $ that generates the topology of $ \Prob (X) $ (recall
\ref{1.10a}) we get a subsequence such that $ \cE{ \phi \circ f }{
\E_{n_k} } \to \cE{ \phi \circ f }{ \E } $ almost everywhere for all $
\phi $ simultaneously; thus $ P_{f|\E_{n_k}} \to P_{f|\E} $ almost
everywhere.
\hfill $ \Box_{\text{claim}} $

\smallskip

\textsc{c. claim.}
If $ \phi : X \to [0,1] $ is upper semicontinuous, then the
function\footnote{%
 Here $ \Ex \phi(f) $ means $ \int_\Om (\phi\circ f) \, dP $.}
\[
L_0 (\Om,\F,P; X) \ni f \mapsto \Ex \, \phi(f) \in \R
\]
is upper semicontinuous.

\textit{Proof of Claim.}
Let $ f, f_1, f_2, \dots \in L_0 (\Om,\F,P; X) $, $ f_n \to f $; we
have to prove that $ \Ex \phi(f) \ge \limsup \Ex \phi(f_n) $. We may
assume that $ f_n \to f $ almost everywhere (not only in probability),
due to the argument of subsequences. Semicontinuity of $ \phi $ gives
$ \phi(f) \ge \limsup \phi(f_n) $ almost everywhere. Therefore $ \Ex
\phi(f) \ge \limsup \Ex \phi(f_n) $ by Fatou's theorem (applied to $ 1
- \phi(f_n) $).
\hfill $ \Box_{\text{claim}} $

\smallskip

\textsc{d. claim.}
The function
\[
\Prob(X) \ni \mu \mapsto \sup_{x\in X} \mu ( \{x\} ) \in [0,1]
\]
is upper semicontinuous.

\textit{Proof of Claim.}
It is easy to check that
\[
\sup_{x\in X} \mu ( \{x\} ) = \inf_{(U_k)} \max_k \mu
 (U_k^{\text{cl}}) \, ;
\]
here the infimum is taken over all finite open coverings $
(U_k)_{k=1,\dots,n} $ of $ X $, and $ U_k^{\text{cl}} $ stands for the
closure of $ U_k $.\footnote{%
 The formula holds also without taking the closure; however, we need
 closed sets here.}
Each $ \mu ( U_k^{\text{cl}} ) $ is upper semicontinuous (in $ \mu $)
by Claim A; therefore the infimum is upper semicontinuous.
\hfill $ \Box_{\text{claim}} $

\smallskip

\begin{sloppypar}
Now we choose some $ f \in L_0 (\Omfav,\Pfav) $ that generates the
whole \mbox{\sif} $ \Ffav $;\footnote{%
 Recalling that $ \Omfav = (0,1) $ we may just take $ f(\om) = \om $.}
a sub-\sif\ $ \E \in \SF^\fav $ is conditionally nonatomic if and only
if almost all conditional distributions $ P_{f|\E} $ are nonatomic,
which may be written as $ \Ex \phi (P_{f|\E}) = 0 $, where $ \phi :
\Prob(\R) \to [0,1] $ is defined by $ \phi (\mu) = \sup_{r\in\R} \mu (
\{r\} ) $. The map $ \SF^\fav \ni \E \mapsto \Ex \phi (P_{f|\E}) $ is
upper semicontinuous, since it is the composition of the map $
\SF^\fav \ni \E \mapsto P_{f|\E} \in L_0 (\Om,\F,P; \Prob(\R) ) $,
continuous by Claim B, and the map $ L_0 (\Om,\F,P; \Prob(\R) ) \ni Z
\mapsto \Ex \phi(Z) $, upper semicontinuous by Claims C and
D. Therefore the set $ \{ \E \in \SF^\fav : \Ex \phi (P_{f|\E}) < \eps
\} $ is open for any $ \eps > 0 $, and so, the set $ \{ \E \in
\SF^\fav : \Ex \phi (P_{f|\E}) = 0 \} $ belongs to the class $ G_\de
$.
\end{sloppypar}
\end{proof}

We may say that comeager many\footnote{%
 Recall \ref{1.3}(c).}
sub-\sif s are conditionally nonatomic, in the following sense.

\begin{corollary}
$ X_1 $ is a dense $ G_\de $, therefore comeager subset of $ \SF^\fav
$.
\end{corollary}

\begin{proof}
By \ref{3.7}, $ X_1 $ is $ G_\de $; also, $ X_1 $ is dense in $
\SF^\fav $, since $ X_1 $ contains all finite sub-\sif s,\footnote{%
 These correspond to \emph{finite} measurable partitions.}
these being dense.
\end{proof}

\begin{theorem}\label{3.10}
(a) $ X_1 $ is a Polish $ \Gfav $-space.

(b) The orbit equivalence relation on $ X_1 $ is smooth.
\end{theorem}

\begin{proof}
(a)
A $ G_\de $-subset of a Polish space is a Polish space, see
\cite[3.11]{Ke} or \cite[2.5.4]{Du}.\footnote{%
 It is about Polish \emph{topological} (not \emph{metric}) spaces.}
We use \ref{3.7} and note that $ X_1 $ is $ \Gfav $-invariant.

(b)
Similarly to the proof of \ref{3.5}(a), for any $ \mu \in \Mfav $
parametrized by $ m_0, m_1, m_2, \dots $ we construct $ \E \in
\SF^\fav $ such that $ (\Omfav,\Pfav) / \E $ is isomorphic to $
(\Omfav,\mu) $, and use an isomorphism between $ (\Omfav,\Pfav) $ and
its square $ (\Om_2,P_2) $, and the conditionally nonatomic sub-\sif\
$ \E_2 = \E \otimes \Zero_{\SF} $. In contrast to that proof, now we
construct $ \E $ in a canonical way, using the fact that $ \Omfav =
(0,1) $. Namely, atoms $ A_k $ of $ \E_k $ are intervals $ A_1 =
(m_0,m_0+m_1)
$, $ A_2 = (m_0+m_1,m_0+m_1+m_2) $, and so on; and $ A_0 = (0,m_0)
$. Thus, $ \E $ depends on $ (m_0,m_1,\dots) $ continuously, and the
set of all these sub-\sif s $ \E $ is closed in $ \SF^\fav $ (in fact,
it is homeomorphic to the compact space $ \Mfav $). The same for the
set of corresponding $ \E_2 $. The latter set, transplanted from $
(\Om_2,P_2) $ to $ (\Omfav,\Pfav) $ by an isomorphism, gives us a
Borel transversal in $ X_1 $ due to \ref{3.1}; recall \ref{1.12}(d)
and \ref{1.15}(a,c).
\end{proof}

So, orbits of $ X_1 $ are in a natural, Borel measurable, one-one
correspondence with points of $ \Mfav $. These are isomorphic types of
probability spaces. Especially, nonatomic probability spaces are
described by a single orbit, consisting of nonatomic $ \E \in X_1
$.

\subsection{Independence and products}

The set $ \SF^\fav $ of sub-\sif s is a lattice; for any $ \E, \F \in
\SF^\fav $ there exist the least sub-\sif\ $ \E \vee \F $ containing $
\E $ and $ \F $, and the greatest sub-\sif\ $ \E \wedge \F $ contained
in $ \E $ and $ \F $. In fact, $ \E \wedge \F $ is just $ \E \cap \F
$; in contrast, $ \E \vee \F $ is generated by $ \E \cup \F $.

If $ \E, \F $ are independent,\footnote{%
 As defined by \ref{2.7}.}
then $ \E \vee \F $ may be called the product of $ \E, \F $ and
denoted also by $ \E \times \F $.\footnote{%
 Or maybe $ \E \otimes \F $.}
(We do not define $ \E \times \F $ when $ \E, \F $ are dependent.) The
same for any finite or countable family of sub-\sif s.

Multiplication of sub-\sif s is closely related to multiplication of
probability spaces, which is well-known (see \cite[Prop.~4]{BEKSY}).

\begin{note}\label{4.1}
(a) Let $ (\Om_1,P_1) $, $ (\Om_2,P_2) $, $ (\Om_3,P_3) $ be three
probability spaces and $ (\Om,P) $ their product. Then $ \E_1 \times
\E_2 = \E_{12} $, where sub-\sif s $ \E_1, \E_2, \E_{12} $ on $
(\Om,P) $ correspond to the first factor, the second factor, and to
both factors, respectively.

(b) Let $ (\Om,P) $ be a probability space, $ \E_1, \E_2, \E_{12} $
sub-\sif s on $ (\Om,P) $ such that $ \E_1 \times \E_2 = \E_{12}
$. Then the quotient space $ (\Om,P) / \E_{12} $ is naturally
isomorphic to the product of two quotient spaces $ ((\Om,P)/\E_1)
\times ((\Om,P)/\E_2) $.
\end{note}

\begin{lemma}\label{4.2}
For every $ \E \in \SF^\fav $ the following three conditions are
equivalent.

(a) $ \E \in X_1 $;

(b) there exists a nonatomic $ \F \in \SF^\fav $ such that $ \E \times
 \F = \One_\SF $;\footnote{%
 That is, $ \E $ and $ \F $ are independent, and $ \E \vee \F $ is the
 whole $ \One_\SF = \Ffav $.}

(c) the morphism $ (\Omfav,\Pfav) \to (\Omfav,\Pfav) / \E $ is
isomorphic to the projection $ (\Omfav\times\Omfav, \mu_1\otimes\Pfav)
\to (\Omfav,\mu_1) $ for some $ \mu_1 \in \Mfav $.
\end{lemma}

\begin{proof}
(a) $ \Longleftrightarrow $ (c) by Theorem \ref{3.2}; and (b) $
\Longleftrightarrow $ (c) by Note \ref{4.1}.
\end{proof}

\begin{lemma}\label{4.3}
For every $ \E_1, \E_2, \dots \in X_1 $ there exist $ \E \in X_1 $ and
$ g_1, g_2, \dots \in \Gfav $ such that
\[
(g_1\cdot\E_1) \times (g_2\cdot\E_2) \times \dots = \E \, .
\]
\end{lemma}

\begin{proof}
The countable product of probability spaces
\[
(\Om_\infty,P_\infty) = (\Omfav,\Pfav) \times (\Omfav,\Pfav) \times
 \dots
\]
is again a probability space isomorphic to $ (\Omfav,\Pfav)
$. Consider sub-\sif s $ \ti\E_1 = \E_1 \otimes \Zero_\SF \otimes
\Zero_\SF \otimes \dots $, $ \ti\E_2 = \Zero_\SF \otimes \E_2 \otimes
\Zero_\SF \otimes \dots $, and so on. They are independent, and $
\ti\E = \ti\E_1 \times \ti\E_2 \times \dots $ is conditionally
nonatomic (since $ \ti\E_1 $ is; use \ref{4.2}(b)). We transplant $
\ti\E_n $ from $ (\Om_\infty, P_\infty) $ to $ (\Omfav,\Pfav) $ by
some isomorphism and apply \ref{3.5}(b), getting $ g_n \cdot \E_n $, $
\E $.
\end{proof}

\begin{note}
The orbit of $ \E $ is uniquely determined by orbits of $ \E_1, \E_2,
\dots $\footnote{%
 Hint: use \ref{3.5}(b).}
Therefore, multiplication of sub-\sif s induces on the set of orbits,
$ X_1 / \Gfav $, an associative operation that takes any finite or
countable number of operands. The corresponding operation on $ \Mfav $
is easy to describe explicitly.
\end{note}

\begin{lemma}\label{4.4}
Let $ \E_1, \E_2, \dots \in \SF^\fav $ be independent. Consider the
countable product
\[
\{0,1\}^\infty = \{0,1\} \times \{0,1\} \times \dots
\]
of two-point topological spaces. For every $ i = (i_1,i_2,\dots) \in
\{0,1\}^\infty $ consider the sub-\sif\
\[
\E_i = \E_1^{i_1} \times \E_2^{i_2} \times \dots \, ,
\]
where $ \E_n^1 = \E_n $, $ \E_n^0 = \Zero_\SF $. Then the map
\[
\{0,1\}^\infty \ni i \mapsto \E_i \in \SF^\fav
\]
is continuous.
\end{lemma}

\begin{proof}
Due to \eqref{2a3} we have to prove continuity of the function
\[
\{0,1\}^\infty \ni i \mapsto \cE{ f }{ \E_i } \in L_2 (\Omfav,\Pfav)
\]
for every $ f \in L_2 $. We may assume that $ f $ is measurable
w.r.t.\ $ \E_1 \times \E_2 \times \dots \, $ Then $ f $ is a sum of
functions of the form $ f_{k_1} \dots f_{k_n} $ where $ k_1 < \dots <
k_n $, and each $ f_{k_j} $ is measurable w.r.t.\ $ \E_{k_j} $, and $
\Ex f_{k_j} = 0 $ for all $ j $. However, $ \cE{ f_{k_1} \dots f_{k_n}
}{ \E_i } = (i_{k_1} \dots i_{k_n}) f_{k_1} \dots f_{k_n} $, which
evidently is continuous in $ i $.
\end{proof}

\begin{corollary}\label{3.17a}
Let $ \E_1, \E_2, \dots \in \SF^\fav $ be independent. Then $ \E_n \to
\Zero_\SF $ and $ \E_1 \times \E_n \to \E_1 $ for $ n \to \infty $.
\end{corollary}

\begin{lemma}\label{3.17b}
For every $ \E \in X_1 $ there exist $ g_1, g_2, \dots \in \Gfav $
such that
\[
g_k \cdot \E \to \Zero_\SF \quad \text{for } k \to \infty \, .
\]
\end{lemma}

\begin{proof}
Lemma \ref{4.3} gives us $ g_k $ such that $ g_k \cdot \E $ are
independent. Corollary \ref{3.17a} ensures that $ g_k \cdot \E \to
\Zero_\SF $.
\end{proof}

\begin{corollary}
The closure of every orbit in $ X_1 $ contains the trivial \sif\ $
\Zero_\SF $.
\end{corollary}

\subsection{Ergodicity}

\begin{lemma}\label{3.18}
$ X_1 $ is ergodic (as defined by \ref{1.153}).
\end{lemma}

\begin{proof}
Due to \ref{1.154} it suffices to prove that $ X_1 $ contains a dense
orbit. We'll see that \emph{the orbit of nonatomic sub-\sif s is
dense.} Let $ \E_0 \in X_1 $ be nonatomic, and $ \E \in X_1 $ be
arbitrary. Lemma \ref{4.3} gives us $ g_1, g_2, \dots \in \Gfav $ such
that $ \E, g_1 \cdot \E_0, g_2 \cdot \E_0, \dots $ are
independent. Corollary \ref{3.17a} shows that $ \E \times g_k \cdot
\E_0 \to \E $ for $ k \to \infty $. However, $ g_k \cdot \E_0 $ is
nonatomic, therefore $ \E \times g_k \cdot \E_0 $ is nonatomic. Being
also conditionally nonatomic, it belongs to the orbit of nonatomic
sub-\sif s. Thus, $ \E $ belongs to the closure of that orbit.
\end{proof}

Combining \ref{3.18}, \ref{3.10}(b) and \ref{1.155} we see that $ X_1
$ contains a comeager orbit. It is easy to guess that it is the orbit
of nonatomic sub-\sif s. The next result confirms the guess.

\begin{proposition}\label{3.21}
The orbit of nonatomic sub-\sif s is a dense $ G_\de $, therefore
comeager subset of $ X_1 $.
\end{proposition}

\begin{proof}
The orbit is dense (as was shown in the proof of \ref{3.18}); we have
to prove that it is $ G_\de $. Similarly to the proof of \ref{3.7}
we'll find an upper semicontinuous function vanishing exactly on that
orbit.

For a given $ f \in L_2 (\Omfav,\Pfav) $ consider the (unconditional)
distribution $ P_{\scE f \E} $ of the conditional expectation $ \cE f
\E $. The function
\[
\SF^\fav \ni \E \mapsto P_{\scE f \E} \in \Prob(\R)
\]
is continuous (recall \eqref{2a3} and \ref{1.10a}). On the other hand,
the function
\[
\Prob (\R) \ni \mu \mapsto \phi(\mu) = \sup_{r\in\R} \mu (\{r\}) \in
 [0,1]
\]
is upper semicontinuous (recall Claim \ref{3.7}.D). Hence the function
\[
\SF^\fav \ni \E \mapsto \phi \( P_{\scE f \E} \) \in [0,1]
\]
is upper semicontinuous. Therefore, its infimum over all $ f \in L_2
(\Omfav,\Pfav) $ is upper semicontinuous. It is easy to see that the
infimum vanishes if and only if $ \E $ is nonatomic.
\end{proof}

\section{Isomorphism of filtrations\\
 as an orbit equivalence relation}
\subsection{The construction}

Recall that filtrations are defined by \ref{2.3}, and their
isomorphisms --- by \ref{2.5}(a).

Lemma \ref{3.5} has a counterpart for filtrations.
Denote by $ X $ the set of all filtrations $ (\E_t) $ on $
(\Omfav,\Pfav) $ such that $ \E_\infty \in X_1 $.

\begin{lemma}\label{3.11}
(a) Every filtration on every probability space is isomorphic to some
filtration belonging to $ X $.

(b) Two filtrations belonging to $ X $ are isomorphic if and only if
some automorphism $ g \in \Gfav $ sends one of them to the other.
\end{lemma}

\begin{proof}
(a) Lemma \ref{3.5}(a) gives us $ \E \in X_1 $ and an isomorphism $
\al $ between the given probability space and $ (\Omfav,\Pfav) / \E
$. The $ \al $-image of the given filtration is an isomorphic filtration
belonging to $ X $.

(b) If $ g $ sends one filtration to the other then they are evidently
isomorphic. On the other hand, let two filtrations $ (\E_t^{(1)}),
(\E_t^{(2)}) \in X $ be isomorphic. Their isomorphism $ \al $ is an
isomorphism between quotient spaces $ (\Omfav,\Pfav) / \E_\infty^{(1)}
$ and $ (\Omfav,\Pfav) / \E_\infty^{(2)} $ sending each $ \E_t^{(1)} $
to $ \E_t^{(2)} $. Note \ref{3.6} gives us $ g \in \Gfav $ that does
the same.
\end{proof}

\begin{note}\label{3.12}
In addition to Item (b) of Lemma \ref{3.11}:

Every isomorphism between two filtrations belonging to $ X $ is
induced by some (at least one) $ g \in \Gfav $.
\end{note}

(A proof is implicitly contained in the proof of Lemma \ref{3.11}.)

Dealing with continuous-time filtrations it is natural to require
continuity (recall \ref{2.3}(b)). On the other hand, discrete-time
filtrations should not be forgotten. It is convenient to deviate from
Definition \ref{2.3} as follows.

\begin{definition}
Let $ T $ be a compact subset of $ [0,+\infty] $.\footnote{%
 Or just a linearly ordered compact topological space, embeddable to $
 [0,+\infty] $ (or equivalently, to $ [0,1] $).}
We define $ X_T $ as the set of all families $ (\F_t)_{t\in T} $ of
 sub-\sif s $ \F_t \in X_1 $ such that
\begin{myitemize}
\item
$ \F_s \subset \F_t $ whenever $ s \le t $;
\item
the map $ T \ni t \mapsto \F_t \in X_1 $ is continuous;
\item
the least sub-\sif\ (corresponding to the least element of $ T $) is
the trivial \sif\ $ \Zero_{\SF} $ (consisting of sets of probability $
0 $ or $ 1 $ only).
\end{myitemize}
\end{definition}

The set $ X_T $ is a subset of the Polish $ \Gfav $-space $
C(T,X_1) $ (recall \ref{1.28}). It is easy to see that $ X_T $ is a closed
invariant set. Therefore
\[
X_T \text{ is a Polish $ \Gfav $-space.}
\]
Lemma \ref{3.11} shows that orbits of $ X_T $ are in a natural
one-one correspondence with isomorphic types of continuous filtrations
on $ T $ (starting with the trivial \sif).

Immersed filtrations and morphisms of filtrations were defined by
\ref{2.4}, \ref{2.5} in the framework of 2b, but their counterparts
for $ X_T $ are evident.

The trivial filtration $ \Zero_{X_T} = ( \Zero_\SF )_{t\in T} $ is
immersed into every filtration. Its orbit is a single point.

\begin{lemma}\label{3.14}
For every filtrations $ x,y \in X_T $, the following two conditions
are equivalent.

(1) There exists a morphism from $ x $ to $ y $.

(2) There exists $ g \in \Gfav $ such that $ g \cdot y $ is immersed
    into $ x $.
\end{lemma}

\begin{proof}
Evidently, (2) implies (1). Assume (1). A morphism from $ x $ to $ y $
is, by definition, an isomorphism between $ y $ and some filtration $
z $ immersed into $ x $; and $ z $ necessarily belongs to $ X_T
$. Note \ref{3.12} gives us $ g \in \Gfav $ such that $ g \cdot y = z
$.
\end{proof}

\begin{lemma}\label{3.15}
Pairs $ (x,y) $ of filtrations $ x,y \in X_T $ such that $ x $ is
immersed into $ y $ are a closed subset of $ X_T \times X_T $.
\end{lemma}

\begin{proof}
As was said after Definition \ref{2.4}, immersion can be expressed in
operator form as $ \F_s \E_t = \E_s $. The operator $ \E_t $ of
conditional expectation depends continuously on the corresponding
sub-\sif\ (recall \eqref{2a3}), therefore, on the filtration $ x
$. The strong operator topology is meant. The product of two operators
is jointly continuous (in these operators), as far as all operators
are of norm $ \le 1 $.
\end{proof}

\begin{lemma}
Pairs $ (x,y) $ of filtrations $ x,y \in X_T $ satisfying equivalent
conditions of Lemma \ref{3.14} are an analytic subset of $ X_T
\times X_T $.
\end{lemma}

\begin{proof}
We take the closed (by \ref{3.15}) set of pairs $ (x,y) $ such that $
y $ is immersed into $ x $, multiply it by $ \Gfav $ and apply the
continuous map $ (g,x,y) \mapsto (x, {g^{-1}\cdot y}) $.
\end{proof}

Especially, we may take $ T = [0,+\infty] $ and consider all
continuous filtrations immersed into Brownian filtrations.\footnote{%
 These are filtrations isomorphic to the natural filtration of the
 Brownian motion, see 2c2.}
The set of such filtrations is analytic, therefore it has the Baire
property and is universally measurable. Similarly we may consider
filtrations that contain immersed Brownian filtrations.

\subsection{Independence and products}

Each filtration $ x \in X_T $ has its maximal sub-\sif\ (corresponding
to the maximal point of $ T $). Filtrations are called independent, if
their maximal sub-\sif s are independent.

For two independent filtrations $ x,y \in X_T $, $ x=(\E_t)_{t\in T}
$, $ y=(\F_t)_{t\in T} $, we define their product $ x \times y $ as
the filtration $ (\E_t\times\F_t)_{t\in T} $; it is again a continuous
filtration;\footnote{%
 Hint: check continuity in $ t $ of $ \cE{ f }{ \E_t\times\F_t } $ for
 the special case when $ f $ is the product of an $ \E_\infty
 $-measurable function and an $ \F_\infty $-measurable function.}
it belongs to $ X_T $ provided that its maximal sub-\sif\ belongs to $
X_1 $. Both $ x $ and $ y $ are immersed into $ x \times y $. The same
for the product of any finite or countable number of filtrations.

\begin{lemma}\label{4.7}
For every $ x_1, x_2, \dots \in X_T $ there exist $ x \in X_T $ and $
g_1, g_2, \dots \in \Gfav $ such that
\[
(g_1 \cdot x_1) \times (g_2 \cdot x_2) \times \dots = x \, .
\]
\end{lemma}

\begin{proof}
Apply Lemma \ref{4.3} to maximal sub-\sif s of given filtrations.
\end{proof}

\begin{note}
Each $ g_k \cdot x_k $ is immersed into $ x $.
\end{note}

\begin{note}
The orbit of $ x $ is uniquely determined by orbits of $ x_1, x_2,
\dots $\footnote{%
 Hint: use \ref{3.11}(b).}
Therefore, multiplication of filtrations induces on the set of orbits,
$ X_T / \Gfav $, an associative operation that takes any finite or
countable number of operands.
\end{note}

\begin{lemma}
For every $ x \in X_T $ there exist $ g_1, g_2, \dots \in \Gfav $ such
that
\[
g_k \cdot x \to \Zero_{X_T} \quad \text{for } k \to \infty \, .
\]
\end{lemma}

\begin{proof}
Apply Lemma \ref{3.17b} to maximal \sif s of given filtrations.
\end{proof}

\begin{corollary}
The closure of every orbit in $ X_T $ contains the trivial filtration
$ \Zero_{X_T} $.
\end{corollary}

\begin{lemma}
For every $ x,y \in X_T $ there exist $ x_1, x_2, \dots \in X_T $ and
$ g_1, g_2, \dots \in \Gfav $ such that $ x_k \to x $, and $ g_k \cdot
y $ is immersed into $ x_k $ for each $ k $.
\end{lemma}

\begin{proof}
Lemma \ref{4.7} gives us $ g_k $ such that $ x, g_1 \cdot y, g_2 \cdot
y, \dots $ are independent. Consider $ x_k = x \times (g_k \cdot y)
$. Corollary \ref{3.17a} ensures that $ x_k \to x $.
\end{proof}

\begin{corollary}\label{4.9}
For every $ x \in X_T $ the set
\[
\{ y \in X_T : \exists g \in \Gfav \; \( \, g \cdot x \text{ is
 immersed into $ y $} \, \) \}
\]
is dense in $ X_T $.
\end{corollary}

For example, filtrations containing immersed Brownian filtrations are
a dense subset of $ X_{[0,+\infty]} $.

\subsection{Ergodicity}

Consider first the case $ T = \{ 0,1,2 \} $, denoting $ X_T = X_{
\{0,1,2\} } $ simply by $ X_2 $. An element of $ X_2 $ may be thought
of as a triple $ ( \E_0=\Zero_{\SF}, \E_1, \E_2 ) $ or just a pair $
(\E_1,\E_2) $ of sub-\sif s $ \E_1,\E_2 \in \SF^\fav $ such that $
\E_1 \subset \E_2 $ and $ \E_2 \in X_1 $. The latter is equivalent to
\[
\E_2 \times \E_{2,\infty} = \One_{\SF} \quad \text{for some nonatomic
 $ \E_{2,\infty} \in \SF^\fav $}
\]
(recall \ref{4.2}). It may happen (or not) that
\[
\E_1 \times \E_{12} = \E_2 \quad \text{for some nonatomic $ \E_{12}
 \in \SF^\fav $} \, ,
\]
in which case we say that $ \E_2 | \E_1 $ is nonatomic.\footnote{%
 Though, $ \E_2 | \E_1 $ itself is undefined.}
It may also happen (or not) that $ \E_1 $ is nonatomic.

\begin{lemma}\label{4.12}
The set $ \{ (\E_1,\E_2) \in X_2 : \E_2 | \E_1 \text{ is nonatomic }
\} $ is a $ G_\de $-subset of $ X_2 $.
\end{lemma}

\begin{proof}
First, a general claim strengthening Claim \ref{3.7}.B.

\smallskip

\textsc{a. claim.} The map
\[
L_0 (\Om,\F,P; X) \times \SF(\Om,\F,P) \ni (f,\E) \mapsto P_{f|\E} \in
 L_0 (\Om,\F,P; \Prob(X) )
\]
is continuous.

\textit{Proof of Claim.}
Let $ f_n \to f $ in $ L_0 (\Om,\F,P; X) $, $ \E_n \to \E $ in $
\SF(\Om,\F,P) $, and $ \phi : X \to \R $ a bounded continuous
function; it suffices to find a subsequence $ (f_{n_k},\E_{n_k}) $
such that $ \cE{ \phi \circ f_{n_k} }{ \E_{n_k} } \to \cE{ \phi \circ
f }{ \E } $ almost everywhere (after that the proof is finished
similarly to the proof of Claim \ref{3.7}.B). We choose $ n_k $
satisfying two conditions: $ \cE{ \phi\circ f }{ \E_{n_k} } \to \cE{
\phi\circ f }{ \E } $ almost everywhere, and $ \sum_k \Ex \, |
\phi\circ f_{n_k} - \phi\circ f | < \infty $. Then
\begin{multline*}
\sum_k \Ex \,\big|\, \cE{ \phi\circ f_{n_k} }{ \E_{n_k} } - \cE{ \phi\circ
 f }{ \E_{n_k} } \big| = \\
\sum_k \Ex \,\big|\, \cE{ \phi\circ f_{n_k} -
 \phi\circ f }{ \E_{n_k} } \big| \le \sum_k \Ex \, | \phi\circ f_{n_k} -
 \phi\circ f | < \infty \, ,
\end{multline*}
\begin{sloppypar}\noindent
thus $ \cE{ \phi\circ f_{n_k} }{ \E_{n_k} } - \cE{ \phi\circ f }{
\E_{n_k} } \to 0 $ almost everywhere. So, $ { \cE{ \phi \circ f_{n_k} }{
\E_{n_k} } } - \cE{ \phi \circ f }{ \E } = \cE{ \phi \circ f_{n_k} }{
\E_{n_k} } - \cE{ \phi \circ f }{ \E_{n_k} } + \cE{ \phi \circ f }{
\E_{n_k} } - \cE{ \phi \circ f }{ \E } \to 0 $ almost everywhere.
\hfill $ \Box_{\text{claim}} $
\end{sloppypar}

\smallskip

The needed nonatomicity may be written as $ \psi ( \E_2 | \E_1 ) = 0
$, where
\[
\psi ( \E_2 | \E_1 ) = \inf_f \, \Ex \phi \( P_{\scE{ f }{ \E_2 } |
\E_1 } \)
\]
and $ \phi $ is the same as in the proof of \ref{3.7} (and
\ref{3.21}); it remains to prove that $ \psi $ is an upper
semicontinuous function on $ X_2 $.

We know from \eqref{2a3} that the map $ \SF \ni \E \mapsto \cE f \E \in
L_2 $ is continuous for every $ f \in L_2 $. Therefore the map $ X_2
\ni (\E_1,\E_2) \mapsto \cE{ f }{ \E_2 } \in L_2 $ is continuous. On
the other hand, Claim A shows that the map
\[
L_2 \times X_2 \ni \( g, (\E_1,\E_2) \) \mapsto P_{g|\E_1} \in L_0
 (\Om,\F,P; \Prob(\R) )
\]
is continuous. Combining the two facts we see that the map
\[
X_2 \ni (\E_1,\E_2) \mapsto P_{\scE{f}{\E_2} | \E_1} \in L_0
 (\Om,\F,P; \Prob(\R) ) 
\]
is continuous. Using Claims \ref{3.7}.C,D we see that the map $ X_2
\ni (\E_1,\E_2) \mapsto \Ex \phi \( P_{\scE{f}{\E_2} | \E_1} \) $ is
upper semicontinuous for every $ f \in L_2 $.
\end{proof}

\begin{definition}
\begin{sloppypar}
A filtration $ (\F_t)_{t\in T} \in X_T $ is called \emph{conditionally
nonatomic,} if $ \F_t | \F_s $ is nonatomic whenever $ s < t $, $ s
\in T $, $ t \in T $. The set of all conditionally nonatomic
filtrations belonging to $ X_T $ is denoted by $ \Xcna_T $.
\end{sloppypar}
\end{definition}

The next result shows that comeager many filtrations are conditionally
nonatomic. It shows also that $ \Xcna_T $ is a Polish $ \Gfav $-space.

\begin{theorem}\label{4.16a}
$ \Xcna_T $ is a dense $ G_\de $-set in $ X_T $.
\end{theorem}

\begin{proof}
We choose a countable subset $ T_0 \subset T $ such that for every $
s,t \in T $ satisfying $ s < t $ there exist $ s_0, t_0 \in T_0 $
satisfying $ s \le s_0 < t_0 \le t $.\footnote{%
 Take a dense subset and add endpoints of intervals that constitute
 the complement of $ T $.}
Nonatomicity of $ \F_{t_0} |
\F_{s_0} $ ensures nonatomicity of $ \F_t | \F_s $. Therefore it
suffices to prove that nonatomicity of $ \F_t | \F_s $ selects a dense
$ G_\de $-set for each pair $ s,t $ separately.

For given $ s,t $ the map $ X_T \ni (\F_r)_{r\in T} \mapsto
(\F_s,\F_t) \in X_2 $ is continuous. Due to Lemma \ref{4.12}, the
considered set is a $ G_\de $-set in $ X_T $. The set is dense, which
follows from \ref{4.9} and existence of a conditionally nonatomic
filtration.\footnote{%
 For example, a Brownian filtration restricted to $ T $.}
\end{proof}

For a finite $ T $, say $ T = \{ 0,1,\dots,n \} $, we denote $ X_T =
X_{ \{ 0,1,\dots,n \} } $ simply by $ X_n $. The simplest infinite $ T
$ is an increasing sequence (and its limit), say $ T = \{ 0,1,2,\dots;
+\infty \} $; here we denote $ X_T $ by $ X_\infty $. The next fact is
well-known.

\begin{proposition}
Each space $ \Xcna_n $ is a single orbit. Also $ \Xcna_\infty $ is a
single orbit.
\end{proposition}

\begin{proof}
We consider $ \Xcna_\infty $ only.
Let $ (\E_n)_n \in \Xcna_\infty $, then $ \E_{n+1} = \E_n \times
\E_{n,n+1} $ for some nonatomic $ \E_{n,n+1} $. Thus $ \E_n = \E_{0,1}
\times \E_{1,2} \times \dots \times \E_{n-1,n} $, which describes $
(\E_n)_n $ uniquely up to isomorphism.
\end{proof}

We return to arbitrary $ T $.

\begin{corollary}\label{4.16}
Let $ x,y \in \Xcna_T $ and $ T_0 \subset T $ be a finite subset. Then
there exists $ g \in \Gfav $ such that $ (g\cdot x) |_{T_0} = y|_{T_0}
$.\footnote{%
 That is, denoting $ x = (\E_t)_{t\in T} $, $ y = (\F_t)_{t\in T} $ we
 have $ g \cdot \E_t = \F_t $ for all $ t \in T_0 $.}
\end{corollary}

\begin{lemma}\label{4.17}
Let $ x \in X_T $ and $ U $ be a neighborhood of $ x $. Then there
exists a finite subset $ T_0 \subset T $ such that
\[
\forall y \in X_T \quad \( \, x|_{T_0} = y|_{T_0} \imply y \in U \, )
 \, .
\]
\end{lemma}

\begin{proof}
Recall the \emph{sandwich argument} stated in \ref{1.14} for $ \bF(X)
$. The space $ \bL(H) $ introduced by \ref{1.10} inherits from $ \bF(H)
$ its topology and (partial) order. Therefore the sandwich argument
holds also for $ \bL(H) $. Further, $ \SF(\Om,\F,P) $ introduced by
\eqref{2a2} inherits its topology and order from $ \bL \( L_2
(\Om,\F,P) \) $. Therefore the sandwich argument holds for $ \SF
$. So, for any $ \E, \E_n, \E'_n, \E''_n \in X_1 $, if $ \E'_n \to \E
$, $ \E''_n \to \E $, and $ \E'_n \subset \E_n \subset \E''_n $ for
all $ n $, then $ \E_n \to \E $.

Recalling that $ \dist (x,y) = \sup_{t\in T} \dist \( x(t), y(t) \) $
(the choice of a compatible metric on $ X_1 $ does not matter), we
choose $ \eps $ such that $ \dist (x,y) \le \eps \impl y \in U
$. Using the sandwich argument and continuity of $ x : T \to X_1 $ we
see that for every $ t \in T $ there exist $ t', t'' \in T $ such that
\begin{gather*}
[t',t''] \text{ is a neighborhood of } t \text{ in } T \, ; \\
\forall \E \in X_1 \quad \( \, x(t') \subset \E \subset x(t'') \imply
 \dist (\E,x(t)) \le \eps/2 \, \) \, .
\end{gather*}
(Of course, the neighborhood $ [t',t''] $ need not be open. We have $
t' \le t \le t'' $; however, do not think that $ t' < t < t'' $; it
may fail, when $ t $ is not an interior point of $ T \subset
[0,+\infty] $.) We choose a finite subcovering and construct $ T_0 $
as the set of endpoints $ t', t'' $ of elements of the subcovering.
\end{proof}

\begin{theorem}\label{4.20}
In the space $ \Xcna_T $, every orbit is dense.
\end{theorem}

\begin{proof}
Let $ x,y \in \Xcna_T $ and $ U \subset \Xcna_T $ be a neighborhood of
$ y $. Lemma \ref{4.17} gives us an appropriate finite $ T_0 \subset T
$. Corollary \ref{4.16} gives us $ g \in \Gfav $ such that $ (g\cdot
x)|_{T_0} = y|_{T_0} $, which implies $ g \cdot x \in U $.
\end{proof}

\begin{theorem}
The space $ X_T $ is ergodic.
\end{theorem}

\begin{proof}
Due to \ref{1.154} it suffices to find a dense orbit. Choose any $ x
\in \Xcna_T $. By \ref{4.20} its orbit is dense in $ \Xcna_T
$. However, $ \Xcna_T $ is dense in $ X_T $ by \ref{4.16a}.
\end{proof}

\stepcounter{footnote}
\footnotetext{%
 A citation, say, \cite[8.4]{Ke} means Item 8.4 (be it a theorem,
 definition or whatever) in the book \cite{Ke}. That is unambiguous
 for \cite{BK}, \cite{Du}, \cite{Ke}, \cite{Sri}. However,
 \cite[2.3]{Hj} could mean either Item 2.3 (in fact, Lemma in
 Sect.~2.1) of \cite{Hj}, or Section 2.3 of \cite{Hj}. I always mean
 an item, unless `Sect.' is indicated explicitly.}

\enlargethispage*{2mm}

\bigskip
\filbreak
{
\small
\begin{sc}
\parindent=0pt\baselineskip=12pt

School of Mathematics, Tel Aviv Univ., Tel Aviv
69978, Israel
\smallskip
\emailwww{tsirel@math.tau.ac.il}
{http://www.math.tau.ac.il/$\sim$tsirel/}
\end{sc}
}
\filbreak

\end{document}